\documentclass[12pt]{article}
\usepackage[english]{babel}
\usepackage{amsmath}
\usepackage{amssymb}
\usepackage{amsthm}
\usepackage{thmtools}
\usepackage{amsfonts}
\usepackage{float}
\usepackage{aliascnt}
\usepackage{fancyhdr}
\usepackage{cleveref}
\usepackage{graphicx}
\usepackage{xcolor}
\usepackage{bbm}
\usepackage{url}
\usepackage{stmaryrd}
\usepackage[margin=1in]{geometry}
\usepackage[utf8]{inputenc}
\usepackage[english]{babel}
\usepackage{tikz}
\definecolor{ku}{RGB}{144,26,30}
\usetikzlibrary{shapes,snakes}
\tikzset{cross/.style={cross out, draw=black, minimum size=2*(#1-\pgflinewidth), inner sep=0pt, outer sep=0pt},
	cross/.default={3pt}}

\makeatletter
\def\thmt@refnamewithcomma #1#2#3,#4,#5\@nil{%
  \@xa\def\csname\thmt@envname #1utorefname\endcsname{#3}%
  \ifcsname #2refname\endcsname
    \csname #2refname\expandafter\endcsname\expandafter{\thmt@envname}{#3}{#4}%
  \fi
}
\newcommand{\indep}{\rotatebox[origin=c]{90}{$\models$}}
\makeatother

\usepackage{etoolbox}
\author{Mads Bonde Raad}
\title{Renewal Time Points for Hawkes Processes}
\newcommand{\N}{\mathbb{N}}
\newcommand{\con}{\rightarrow}
\newcommand{\conx}[1]{\stackrel{#1}{\rightarrow}}
\newcommand{\R}{\mathbb{R}}

\newcommand{\define}[4]{\expandafter#1\csname#3#4\endcsname{#2{#4}}}

\forcsvlist{\define{\newcommand}{\mathcal}{c}}{A,B,C,D,E,F,G,H,I,J,K,L,P,M,N,U,S,T,W,R,V,X,Y}
\forcsvlist{\define{\newcommand}{\mathbb}{q}}{A,B,C,D,E,F,G,H,I,J,K,L,P,Q,M,Z,T,X,Y}
\newcommand{\ql}{\mathbbm{1}}

\newcommand{\lp}{\left( }
\newcommand{\rp}{\right) }
\newcommand{\lb}{\left\lbrace}
\newcommand{\rb}{\right\rbrace}
\newcommand{\lf}{\left[}
\newcommand{\rf}{\right]}
\newcommand{\lff}{\lfloor}
\newcommand{\rff}{\rfloor}
\newcommand{\lc}{\left\lceil}
\newcommand{\rc}{\right\rceil}
\newcommand{\lv}{\left\vert}
\newcommand{\rv}{\right\vert}
\newcommand{\lV}{\left\Vert}
\newcommand{\rV}{\right\Vert}
\declaretheoremstyle[
  spaceabove=\topsep,
  spacebelow=\topsep,
  headfont=\normalfont\bfseries,
  notefont=\bfseries, notebraces={(}{)},
  bodyfont=\normalfont\itshape,
 postheadspace=\newline,
 qed=$\circ$
]{bombom}
\declaretheoremstyle[
spaceabove=\topsep,
spacebelow=\topsep,
headfont=\normalfont\bfseries,
notefont=\bfseries, notebraces={(}{)},
bodyfont=\normalfont,
postheadspace=\newline,
qed=$\circ$
]{exstyle}

\declaretheorem[style=bombom,title=Theorem,numberwithin=section,refname={theorem}]{s}

\declaretheorem[style=bombom,name=Remark,sibling=s,refname={remark}]{rem}
\declaretheorem[style=bombom,name=Definition,sibling=s,refname={definition}]{defi}
\declaretheorem[style=bombom,name=Corollary,sibling=s,refname={corollary}]{koro}
\declaretheorem[style=exstyle,name=Example,sibling=s,refname={example}]{eks}
\declaretheorem[style=exstyle,name=Assumption,refname={assumption}]{as}
\declaretheorem[style=exstyle,name=Setup,numbered=no,refname={setup}]{set}
\declaretheorem[style=bombom,name=Proposition,sibling=s,refname={proposition}]{prop}

\declaretheoremstyle[
  spaceabove=\topsep,
  spacebelow=30pt,
  headfont=\normalfont\bfseries,
  notefont=\mdseries, notebraces={(}{)},
  bodyfont=\normalfont,
 postheadspace=\newline,
 numbered=no,
  qed=\qedsymbol,
  name=Proof
]{mythmstyle}
\declaretheorem[style=mythmstyle]{p}
\begin{document}
\maketitle
\begin{abstract}
	In the last decade Hawkes processes have received much
        attention as models for functional connectivity in neural
        spiking networks and other dynamical systems with a cascade behavior. In this paper we establish a renewal approach for analyzing this process. We consider the ordinary nonlinear Hawkes process as well as the more recently described age dependent Hawkes process. We construct renewal-times and establish moment results for these. This gives rise to study the Hawkes process as a Markov chain. As an application, we prove asymptotic results such as a functional CLT and a time-average CLT.
\end{abstract}


\section{Introduction}
\subsection{Hawkes Processes}
 The Hawkes process has received much attention in the last decade for modeling events that exhibit self excitation - or inhibition. There are many
examples of phenomena of interest which exhibit such behavior, including finance that propagate through a market giving rise to volatility clustering
observations \cite{Hawkes-In-Finance},
interactions on social media \cite{Social-Media}, and pattern dependencies in DNA \cite{pat}.

The primary application in mind for this article is neuroscience where the Hawkes process may be used to model spike trains for one or more neurons.
When neurons send an electric signal,
the so-called {\em action potential} or {\it spike}, they excite or inhibit recipient neurons in
the network (the {\em post-synaptic} neurons).  Jumps of the $i$th
unit of the Hawkes
process are then identified with the spike times of the $i$th
neuron. Moreover, the biological process imposes a strong
self-inhibition on a neuron that has just emitted a spike. This period
of about $2 ms$ is called the \textit{absolute refractory period}, and
in this phase it is virtually impossible for a neuron to spike
again. The neuron then gradually regains its ability to spike in the
longer \textit{relative refractory period}. It was proposed in
\cite{chevallier} to model absolute and relative refractory periods in
neuronal spike trains by {\it age dependent Hawkes processes}, where
the age of a unit is defined as the time passed since the last time it
jumped, and thus, it resets to zero at each jump time.

In the present paper we consider a one-dimensional Hawkes process. It may be an ordinary nonlinear type as in \cite{bm}, or an age dependent one as in \cite{chevallier,MES}. The Hawkes process may be described as a random counting measure $Z$ on $\R_{+}$ with an associated intensity process $\lambda$. The intensity is colloquially speaking the conditional probability of a jump of $Z$, given the history ${\mathcal F}_{t}$ of $Z$,
$$  \lambda_t  dt \approx P ( Z\mbox{ has a jump in } \left( t ,
  t + dt \right]\vert {\mathcal F}_{t}).$$
The age dependent Hawkes process is characterized by having an intensity which is a function of a weighted average of the time since past jumps, i.e.:
$$ \lambda_{t}=\psi\lp \int_{0}^{t-}h\lp t-s\rp dZ_{s},A_{t} \rp  $$
where $\psi$ is a function  which is  Lipschitz in the first coordinate, and the so called \textit{age process} $A$  is the time since the last jump of $Z$. In general we denote this process as the "age dependent Hawkes process" (ADHP). It was introduced by \cite{chevallier} and its stability properties were dicussed in \cite{MES}. If $\psi$ does not depend on its 2nd coordinate, we obtain the "nonlinear Hawkes process", which we in the article shall denote as the ordinary Hawkes process to distinquish it from the age dependent Hawkes process. Many aspects of this process have been studied by various authors. Stability was discussed in \cite{bm}, a CLT result was discussed in \cite{Zhu1}, and a wide range of mean-field results have been discussed in \cite{dfh}. A multiclass setup was discussed in \cite{SusEva} focusing on mean-field limits and oscillatory behavior. Finally, if we take $\psi=c+Lx_{+}$, we obtain the Linear Hawkes process. This was the process studied first by Alan G. Hawkes in \cite{Hawkes}. It can be represented as a Poisson Branching process where the centre process is a homogeneous Poisson process of intensity $c$, and the offspring processes are inhomogenous Poisson processes of intensity $h$. See also \cite{patroy,DALEY}.

\subsection{Purpose and Results}
 In this paper we discuss stability of Hawkes processes from a renewal perspective. When $h$ is of compact support and $Z$  is an ergodic linear Hawkes process, it will happen infinitely often that $Z\lf t-\text{supp}\lp h\rp ,t\rf=0$, at which point a renewal occurs. It was shown in \cite{graham} that these renewal times have exponential moment under certain regularity assumptions. However, when the weight function $h$ does not have compact support, it is no longer straightforward to find timepoints where the past can be eliminated. In this article we show how to construct such renewal times. The procedure is not unlike the Athreya-Ney technique for Markov Chains in the sense that we wait for some stopping-time $\alpha_{0}$ to occur, which may be interpreted as a minorization criteria. Here we let random variables independent of $Z$ decide whether we obtain a renewal $\alpha_{0}$ at this point, or we jump to a new state of $Z$ by moving time forward to a stopping time $\tau_{1}$. This procedure is repeated, until a renewal has occured after a random number of iterations $\eta$. The renewal time $\alpha_{\eta}$ will be a stopping time w.r.t. the enlarged filtration induced by $Z$ and the independent decisions.\\\\

The renewal approach to discussing stability of Hawkes processes turn out to be beneficial for establishing a number of key results for Hawkes processes. Here we give a brief overview of the results:
\begin{itemize}
	\item It is well known from \cite{bm} that two Hawkes processes driven by the same Poisson random measure with  sufficiently fast decaying initial signals couple eventually. The coupling time is bounded by the renewal time $\alpha_{\eta}$. We use this to formulate moment results for the coupling time in terms of the distribution $h \ dt$. Moreover, $\alpha_{\eta}$ is constructed explicitly so that it can be simulated.

	\item We prove a CLT for processes of the time average type:
	\begin{align}
		t^{-1/2 }\int_0^t H(Z_{\mid [s-D,s]})ds\Rightarrow N\lp\mu,\sigma^{2} \rp
	\end{align}
for appropiate $\mu\in \R,\sigma>0$.
	This was done for the linear Hawkes processes in \cite{graham} assuming compact support of $h$, and for such $h$ our results coincide.

	\item We prove a functional CLT for Hawkes processes. This was done for ordinary Hawkes processes in \cite{Zhu1} with slightly weaker integrability assumptions on $h$ compared to what we impose. However, we do not need positivity of $h$, nor do we need that $h$ itself is decreasing.
\end{itemize}


\section{Notation, Definitions and Core Assumptions}
Throughout this article, we will be working on a background probability space $\lp \Omega, {\mathcal F},P\rp $ and  all random variables are assumed to be defined on this space. A random variable $X$ is said to have $q$'th moment for some $q\geq 0$ if $\qE \lv X^{q}\rv <\infty$, and it is said to have exponential moment if $\qE \exp\lp c X\rp<\infty$ for some $c>0$. Likewise, we say that a function $f:\R\con \R$ have   $q$'th moment, respectively exponential moment, if $\int x^{q} f\lp x\rp dx<\infty,$ respectively $\int \exp\lp cx\rp  f\lp x\rp dx<\infty$.
We recall the basic Stieltjes integration notation. A function of finite variation $f:\R\con \R $ induces a Stieltjes signed measure $\mu_{f}$ satisfying $\mu_{f}\lp \lp a,b\rf\rp = f\lp b\rp-f\lp a\rp$.  We use the notation $\lv df \rv$ for the corresponding variation measure. The \textit{Lebesgue-Stieltjes integral} is defined as
\begin{align*}
\int_{a}^{b} f\lp x\rp dg\lp x\rp = \int \ql \lb \lp a,b\rf
\rb (x) f\lp x\rp d\mu_{g}\lp x\rp,
\end{align*}
see e.g. \cite{Halmos}. If $\nu$ is a measure on $\R^{2}$ we shall also use the following notation for the integral over semi-closed boxes
\begin{align*}
\int_{a}^{b}\int_{c}^{d} f\lp x,y\rp d\nu\lp x,y\rp = \int \ql \lb \lp a,b\rf \rb \lp x\rp \ql \lb \lp c,d\rf \rb \lp y\rp f\lp x,y\rp d\nu\lp x,y\rp.
\end{align*}
Recall that if $\lp A,d\rp $ is a complete separable metric space (c.s.m.s) then the space of locally bounded counting measures $M^c_{A}$ is a c.s.m.s as well, when equipped with an appropiate metric. See chapter 9 and/or appendix 2 in \cite{DALEY} for an overview of properties for this space. For $n\geq 0$ we introduce the shift operator $\theta$ from the space $M^c_{\R\times \R^n}$ onto itself, as the map
\begin{align}
\theta^{t}\nu\lp C\rp=\nu\lp C+te_{1}\rp,\quad  \forall C\in \cB^{n+1},\nu \in M^c_{\R\times \R^n}
\end{align}
 where $e_{1}\in \R \times \R^{n} $ is the first unit vector. For a measure $\nu\in M^c_{\R\times \R^n}$ and $t\in \R$ we also define the increment measure $\nu_{t+}\in M^c_{\R_{+}\times \R^n}$ by
 \begin{align}
 	\nu_{t+}(C) = \theta^t\nu(C),\ \forall C\in\cB_{\R_+\times \R^n}.
 \end{align}
For the sake of clear notation, we shall agree that $\theta^r \nu_{t+}\lp C\rp:=\lp \theta^r\lp \nu_{t+}\rp\rp \lp C\rp.$\\\\
Let $T\subset\R$ (with possible equality) be an interval and let $\nu$ be a locally bounded measure on $T\times \R_+$. The random variable $\Pi:\Omega\to M^c_{T\times \R_{+}}$ is called a Poisson Random Measure (PRM) with mean measure $\nu$ if
\begin{enumerate}
	\item $\Pi(A) \sim Pois(\nu(A)),\ \forall A\in \cB_{T\times \R_+},\ \nu(A) <\infty$,
	\item $A_1, ..., A_m\in \cB_{T\times \R_{+}}$ disjoint $\Rightarrow$ $\Pi(A_1)\indep ... \indep \Pi(A_m)$.
\end{enumerate}
It is assumed that the mean measure is the Lebesgue measure unless otherwise mentioned. Moreover, let $(\cG_t)_{t\in T}$ be a filtration such that $\Pi(A\cap (-\infty,t]\times \R_+)$ is $\cG_t$ measurable for all $A\in \cB_{T\times \R_+}$. We call $\Pi$ a $\cG_t$-PRM, if $\Pi_{t+}$ is a PRM such that $\Pi_{t+}\indep \cG_t$ for all $t\in T$.
\\In the following we introduce the core mathematical objects and assumptions needed to discuss the Age Dependent Hawkes process.
\begin{itemize}
	\item[$\pi\; \vert $] $\pi $ and $\overline{\pi}$ are independent  PRMs on $\R \times \R_{+} $ with Lebesgue intensity measure.
	\item[$\lp \cF_{t}\rp\; \vert $] We assume $\lp \cF_{t}\rp_{t\geq 0} $ is a filtration such that $\pi,\overline{\pi}$ are $ \cF_{t}$-PRMs.

	\item[$h\; \vert $] The weight function $h:\R_{+}\rightarrow \R$  is a locally integrable function.
	\item[$R\; \vert  $]  The initial signal $(R_{t})_{t \geq 0}$ is an  $\mathcal{F}_{0 }\otimes \mathcal{B} $ measurable process on $\R_{+}$ such that\\  $\qE \int_{0}^{t} |R_s| ds<\infty$ for all $t\geq 0$.
	\item[$\psi\; \vert  $]  The rate function   $\psi : \R\times \R_{+}\con \R_{+}$ is a measurable function which is increasing in both variables and satisfying a Lipschitz-like condition : For all $x\leq y \in\R$ and $a\leq b \in\R_+$ it holds that
	\begin{align}\label{eq:psias1}
		 \psi\left( y,b\right) - \psi\left( x,a\right)\leq \begin{cases}
L \lp y - x\rp &a=b \\
L \lp y - x\rp+ \left( c^{pre}_{\psi}+ L x_{+} \right) g\left( a\right) &a<  b
\end{cases}
	\end{align}
	for some constants $L,c^{pre}_{\psi}>0$ and a decreasing function $g$ bounded by $1$.
	\item[$A_{0}\; \vert $] The initial age $A_{0}$ is a $ {\mathcal F}_0-$measurable random variable with support in $\R_{+}$.
\end{itemize}
We observe that \eqref{eq:psias1} implies that $\psi$ is sublinear since
\begin{align}
&\psi\lp y,b\rp\leq \psi\lp 0,b\rp\pm\psi\lp 0,0\rp\leq c^{pre}_{\psi}g\lp 0\rp +\psi\lp 0,0\rp \quad \text{if   } y<0,b\in \R_{+}\\
&\psi\lp y,b\rp= \psi\lp y,b\rp\pm\psi\lp 0,0\rp\leq Ly+c^{pre}_{\psi}g\lp 0\rp +\psi\lp 0,0\rp \quad \text{if   } y\geq 0,b\in \R_{+}
\end{align}
so with  $c_{\psi}:=c_{\psi}^{pre}+\psi\lp 0,0\rp $ we have
\begin{align}\label{eq:psias2}
\psi\lp y,b\rp\leq c_{\psi}+Ly_{+}\quad \forall y\in \R,b\in \R_{+}.
\end{align}
We may now define the ADHP. For convenience in proofs we also define the  $D$-delayed ADHP for $D\geq 0$, which is essentially an ADHP where the intensity is killed until time $D$.
\begin{defi}
	The  Age Dependent Hawkes Process (ADHP) $Z^*$ driven by $\pi$, with weight function $h$, rate function $\psi$, initial signal $R^*$ and initial age $A^*_{0}$ is the random measure  $\R_{+}$ satisfying the folowing: For all $a\leq b\in \R_{+},t\in \R_{+}$
 \begin{align}
	Z^*(a,b]&= \int_{a}^{b}\int_{0}^{\infty} \ql\{z\leq \lambda^*_s\} d\pi (s,z)\\
	\lambda^*_{t}&=\psi\lp X^*_{t},A^*_{t}\rp\\
	X^*_{t} &= \int_{0}^{t-}h\lp t - s\rp  dZ^*_{s}+R^*_{t}
	\end{align}
	and $A^*$ is the  càdlàg age process of $Z^*$, given as $t-\sup \lb s<t: Z^*\lf s,s\rf=1  \rb$. We refer to $A^*$ as the \textit{Age} of $Z^*$, $X^*$ as the \textit{Memory} of $Z^*$ and $\lambda^*$ as the \textit{Intensity} of $Z^*$.\\
Let $D\geq 0$. We shall say that $Z^*$ is a $D$-delayed ADHP with weight function $h$, rate function $\psi$, initial signal $R^*$ and initial age $A^*_{0}$ if $Z^*\lf t,t\rf =0$ for $t\in \lf 0,D\rf$ and $Z^*_{D+}$ is an ADHP (driven by $\pi_{D+}$) with parameters $\lp h,\psi\rp$, initial age $A^*_{0}+D$ and signal $t\mapsto R^*_{t+D}.$
\end{defi}
The proof of theorem 0.1 \cite{MES} implies that the ADHP is indeed well-defined.
When $D=0$ we obtain the regular ADHP. If it also holds that $\psi\lp x,a\rp $ does not depend on $a$, then it is the ordinary nonlinear Hawkes process. If moreover $\psi\lp x\rp=\psi\lp 0 \rp +Lx_{+}  $ is linear, then we obtain the linear Hawkes process.

Let now $Z^{*}$ be the ADHP that we wish to obtain a regeneration point for. It is well known that depending on the parameters $h,\psi,$ Hawkes processes can either be in the subcritical regime where $\limsup_{t\con \infty} Z^*\lp 0,t\rf/t<\infty$ or in the supercritical regime where the limit is $\infty$.  To succeed we must ensure that $Z^*$ is  in the subcritical regime. We shall treat two different setups that will ensure this. The first setup assumes that $\int_{\R_{+}} h^+\lp s\rp  ds<L^{-1}.$ For the Linear Hawkes process with
\begin{align}
\psi_{L}\lp x\rp:=c_{\psi}+Lx_{+} \label{linear}
\end{align}
this has the interpretation that each direct child of a parent jump induces $<1$ new child on average. The second setup assumes that the ADHP has a refractory period, i.e. the intensity is bounded for a period after each jump (this includes the case where $\psi$ is uniformly bounded).
\begin{set}[Ordinary Hawkes process]
\label{as:systemo}
We assume that   $\int_{\R_{+}} h_{+}\lp s\rp  ds<L^{-1}$ where $h_{+}\lp t\rp =\max\lp h\lp t\rp ,0\rp $.
\end{set}

\begin{set}[Age dependent Hawkes process]
\label{as:systemad}
There exist $K\geq 0,\delta\in \lb 1/n:n\in \N\rb$ s.t.\footnote{It is merely for mathematical convenience that we restrict $\delta$ to reciprocal integers instead of arbitrary $\delta\in \R_{+}$ in setup (AD).}
\begin{align}
\psi\lp x,a\rp\leq K \text{ for } a\in \lf 0,\delta \rf
,x\in \R.\label{eq:refractory}
\end{align}
 \end{set}
We shall establish a renewal time for each of these setups. While some variables will vary slightly in their definitions for each setup, the approach is similar so the renewal time will be constructed simultaneously.  We shall refer to the two setups above as setup (O) or (AD) respectively.

\begin{eks}
  Consider the rate function given by
	\begin{align}
	\psi(x,a)=l(x)\varphi\lp x,a\rp
	\end{align}
  where $l,\varphi$ are increasing, $l$ is Lipschitz and $\varphi$ is bounded by 1 for all $x\in\R$, $a\in \R_{+}.$  Moreover, we assume $\varphi$ converges to $1$ in the sense that there is a function $g:\R_{+}\con \lf 0,1 \rf$ satisfying $1-\varphi\lp x,a\rp \leq g\lp a\rp $ for all $ x\in \R,a\in \R_{+}$. For $\varphi\equiv 1$ we obtain the ordinary Hawkes process, so for general $\varphi$ we may interpret the ADHP as an ordinary Hawkes process with rate function $l$, but inhibited by its own age process with a factor $\varphi\lp x,a\rp$.  \\\\
To show that it satisfies \eqref{eq:psias1} we take arbitrary $x\leq y$ and $a\neq b$ and obtain
	\begin{align}
	\psi(y,b)-\psi(x,a)&=l\lp y\rp \varphi\lp y,b\rp-l\lp x\rp \varphi\lp x,a\rp\\
&= \lp l\lp y\rp-l\lp x\rp\rp  \varphi\lp y,b\rp+l\lp x\rp \lp \varphi\lp y,b\rp-\varphi\lp x,a\rp \rp
	\end{align}
If $x\geq 0$ then $l\lp x\rp\leq l\lp 0\rp+ L\lp x-0 \rp  $ while $l\lp x\rp\leq l\lp 0\rp$ for $x<0$. We use this, and the fact that $\varphi\lp y,b\rp-\varphi\lp x,a\rp\leq 1-\varphi\lp x,a\rp \leq g\lp a\rp$ to conclude

\begin{align*}
\psi(y,b)-\psi(x,a)&\leq L\lp y-x\rp+\lp l\lp 0\rp +Lx_{+}\rp g\lp a\rp
\end{align*}
which fits into \eqref{eq:psias1}.
 The most principal example of $\varphi$ is the simple $\ql \lb A\leq \delta\rp $ corresponding to a hard refractory period, and in this case $Z^*$ is in setup (AD).  Although this is a rather simple example, it is important due to its application for modelling neural spike-trains.

  As mentioned previously, if $\varphi\equiv 1$ then
	\begin{align}
	\psi(x,a) = l(x)
	\end{align}
and one obtains the ordinary Hawkes process. We may assume that $\lV h_{+}\rV_{\cL^1}<L^{-1}$ in which case the parameters fits under setup O, or we can assume $\psi$ is bounded in which case it fits under setup AD.
\end{eks}

For each setup, we impose two assumptions. The first one restricts the randomness in the initial signal.

\begin{as}
\label{As:alpha0}
There is an a.s. finite $\cF_{t}$-stopping time $\alpha_{0}$ and a deterministic decreasing function $r:\R_{+}\con \R_{+}$ such that for all $t\geq \alpha_{0}$
	\begin{align}
	&\lv \int_{0}^{\alpha_{0} } h\lp t - s\rp  dZ^*_{s}+ R^{*}_{t} \rv  \leq r(t-\alpha_{0}) .\label{OffsetEq}
	\end{align}
\end{as}

The next assumption puts integrability assumptions on $r,h,g$. It will be split in two. One where $h,r$ have power tails, and one where they have exponential tails.
\begin{as}
\label{As:fmap}
Let $\gamma :\R_{+}\con \R_{+}$ be an increasing and right continuous function and define $\overline{h}(t) := \sup_{s\geq t} \lv h(s) \rv.$
 We assume that $\overline{h}\in \cL^1_{loc}$ and either (A) or (B) below holds.  \\[5pt]
\textit{Assumption 2 (A):}\\
There exists $p\geq 0$ s.t.
\begin{itemize}
	\item $  t\mapsto t^{p}r\lp t\rp \in \cL^{1}, t\mapsto t^{p} g\lp t\rp \in \cL^{1}, \text{ and } t\mapsto t^{p+1}\gamma\lp t + 1 \rp \overline{h}\lp t\rp\in \cL^{1}$.
  \item Under setup (O) we assume
\begin{align}
\liminf_{t\con \infty}\dfrac{\gamma\lp t\rp  }{c^{-1}_{h}\lp p+1\rp \ln_{ + }t}>1 \label{criteriaO}
\end{align}
 where $c_{h}= \lV h_{+} \rV_{\cL^{1}}-\ln_{+} \lV h_{+} \rV_{\cL^{1}} -1.$ Under (AD) we merely assume

\begin{align}
\liminf_{t\con \infty}\dfrac{\gamma\lp t\rp  }{\ln_{ + }t}>0.
\end{align}

\end{itemize}
\textit{Assumption 2 (B):}
\begin{itemize}
	\item The functions $r,g$ and $\overline{h}$ have exponential moments.
\item We assume that
$$  \liminf_{t\con \infty} \dfrac{\gamma\lp t\rp}{t}>0.   $$
\end{itemize}
\end{as}

\begin{rem}
\begin{enumerate}
A few remarks on the introduced variables and assumptions:
\item For all results to come in this article we shall implicitly assume \cref{As:alpha0} and \cref{As:fmap} (A), unless otherwise stated. We will state explicitly which setup we work under. \Cref{As:fmap} (B) is clearly stronger than the (A) version for any choice of $p$, and we state explicitly when we work under assumption (B) instead of (A).
\item  The map $\overline{h}$ is the smallest decreasing function dominating $h$. As in \cite{MES}, we put integrability assumptions on $\overline{h}$ which is slightly more restrictive than if they were put on $h$. It turns out to be advantagous to work with a decreasing weight function and the restriction is, at least in the belief of the author, of small consequence for practical applications.
\item If $\int_{0}^\infty t^{p+1}\ln_{+} t\; \overline{h}\lp t\rp dt<\infty$ then the choice $\gamma\lp t\rp = C\ln_{+}\lp t\rp $ for large $C$ satisfies the parts of \cref{As:fmap} (A) relevant to $\overline{h},\gamma$.  Likewise, if $\int_{0}^\infty \exp\lp ct\rp \overline{h}\lp t\rp dt<\infty$ for some $c>0$, then the choice $\gamma\lp t\rp=Ct$ for any $C>0$ satisfies \cref{As:fmap} (B).We allow $\gamma$ to be chosen freely because it may change the speed of computation in an actual simulation of the renewal-times. Recall that since $\gamma$ is right continuous, the generalized inverse
$\gamma^{-1}\lp t\rp:=\inf \lb s\geq 0 : \gamma\lp s\rp\geq t  \rb $ satisfies
\begin{align}
y\leq \gamma\lp t\rp \Leftrightarrow \gamma^{-1}\lp y\rp\leq t.
\end{align}
\end{enumerate}
\end{rem}
We now define some key functions to be used in the construction of a regeneration time, and with \cref{As:fmap} we immediately determine their integrability properties. Define
\begin{align}
f\left(t_{1},t_{2}\right)   &= \lp \gamma\lp 0\rp  +1 + \delta^{ - 1} \rp  \lp \overline{h}\left( t_{1} \right)+ \int_{0}^{t_{2}}\gamma\lp s+1\rp  \overline{h}\left( t_{1}  + s\right) ds\rp +r\lp t_{1} \rp
\end{align}
(with convention $\delta = \infty$ in the (O)-system). For convenience we write $f\lp t\rp $ instead of $f\lp t,\infty\rp $. Define also
\begin{align}
F^{pre}\left( t\right) &=2Lf\left( t\right) + c_{\psi}g\left( t\right)
\\F(t) &= \ql\{t\leq D\} \lp c_{\psi}+Lf\lp t\rp\rp   + \ql\{t>D\}F^{pre}(t)
\end{align}

\begin{prop}\label{prop:fmap}
Consider the maps $t\mapsto f\lp t\rp ,F^{pre}\lp t\rp ,F\lp t\rp$ defined above. Under either setup and assumption 2A these functions have p'th moment, and under assumption 2B these functions have exponential moments.
\end{prop}

%
%


An  important example for which \cref{As:alpha0} is satisfied is the stationary ADHP:

\begin{eks}
\label{eks:inv}
 The classical method of studying stability of Hawkes processes, due to Brémaud \& Massoulié \cite{bm} has been to find a solution
\begin{align*}
Z^I(a,b]&=\int_{a}^{b}\int_{0}^{\infty}  \ql\{z\leq \lambda^I_s\} d\pi (s,z)\\
X^I_{t}&=\int_{-\infty}^{t-} h\lp t-s\rp dZ^I_{s}\\
\lambda^I_s&=\psi\lp X_{t},A_{t}\rp
\end{align*}
such that $\theta^t Z^I =H^I\lp \theta^t\pi \rp$ for some suitable map $H^I:M^c_{\R\times \R_{+}}\con M^c_{\R}$ .  Note that $Z^I$ is an ADHP with $R_{t}^I=\int_{-\infty}^{0-} h\lp t-s\rp dZ^I_{s}.$ See \cite{bm,MES} for criterias of existence. In both of the cited papers, it is proven that when $Z^I$ exists, it is stationary and ergodic, and if $Z^*$ is another Hawkes process driven by $\pi$,  with a signal $R^*$ satisfying $\qE \int_{0}^\infty \lv R^*_{s}\rv ds<\infty$, then $Z^*$ couples with $Z^I$ eventually.
  We shall see in \cref{Prop:AlphaExists} that there is a suitable choice of $\alpha_{0}$,$r$ such that \cref{As:alpha0} is satisfied for $Z^I$. In  \cref{CouplingBound} we prove that the coupling time has $p$'th moment and even exponential moment under assumption 2B.
\end{eks}

\section{Renewal for Hawkes Processes}
The purpose of this section is to develop a renewal time point $\rho$ for a given Hawkes process, which we do in section \ref{renewalsection}. In section \ref{markovsection} we use this to write the Hawkes process as a function of a Markov process.
\subsection{Constructing a Renewal Time Point for a Hawkes Process}\label{renewalsection}
 In this section we are given an ADHP $Z^*$.  The goal is to prove the main result \cref{regenerationTheorem} which gives a random time $\rho$ satisfying that $Z^*_{\rho+}\indep Z_{\vert \lp 0,\rho \rf}$. This is done by introducing a point process $Z$ which regenerates at stopping times $\alpha_{n}$. Then, by using the specific construction of $\alpha_{n}$ and $Z$, we are able to throw a biased coin deciding whether $\alpha_{n}$ should be the renewal time point for $Z^*$ or not.\\

 The first step is to construct $Z$. Given either of the two setups,  we shall simultaneously define $Z$, the sequences of stopping times $\lp \tau_{n}\rp, \lp \alpha_{n}\rp$, and intensities $\underline{\lambda},\overline{\lambda}$ as a system. While the system is defined slightly different for each of the two setups, they are very similar. Only the $\alpha$'s differs in the definition, depending on whether we discuss the (O) system or the (AD) system.
 Recall the split PRMs $\pi^{\uparrow},\pi^{\downarrow}$ defined in the first section of the appendix. We also note that we make use of the convention $\inf\{\emptyset\} = \infty$.\\

The system is defined as follows: Fix $D\geq 0$ and define $\overline{\lambda}_{t}=\underline{\lambda}_{t}=0$ for $t\in \lp 0, \alpha_{0} \rf $ and $Z\lp 0,\alpha_{0}\rf=0.$ For $n\in\N$ s.t. $\alpha_{n-1}<\infty$  we define $Z$ as the $D$-delayed ADHP driven by $\pi$, with parameters $h,\psi$ and initial conditions
$A_{\alpha_{n - 1} + } = 0$, $R_{t} = -f\left( t-\alpha_{n - 1}\right)$.
Let $\lambda$ be its intensity. We set $\underline{\lambda}_{t},\overline{\lambda}_{t}=0$ if $t\in \lp \tau_{n},\alpha_{n}\rf$  and
\begin{align}
\underline{\lambda}_{t}
&=\lambda_{t}\\
\overline{\lambda}_{t}
&=\underline{\lambda}_{t} + F(t-\alpha_{n-1})
\end{align}
when  $t\in \lp \alpha_{n - 1},\tau_{n}\rf$. Moreover, we set
\begin{align}
\tau_{n}
&=\inf\left\{ t>\alpha_{n - 1}:\int_{\alpha_{n - 1}}^{t}\int_{0}^{\infty}\ql\left\{ z \in \left( \underline{\lambda}_{s},\overline{\lambda}_{s} \right]\right\} d\pi\left( s,z\right) \geq 1\right\}&\\
&=\inf\left\{ t>\alpha_{n - 1}:\int_{\alpha_{n - 1}}^{t}\int_{0}^{\infty}\ql\left\{ z \leq F\lp t-\alpha_{n - 1}\rp \right\} d\pi^{\uparrow \underline{\lambda},\overline{\lambda}}\left( s,z\right) \geq 1\right\}.
\end{align}
If $\tau_{n}=\infty$ we set $\alpha_{n}=\infty$ in either setup. Otherwise, under setup (AD) we choose
\begin{align}
  \alpha_{n} &= \alpha_{n-1}+\inf\{ i>\lceil \tau_{n}-\alpha_{n-1}\rceil: \theta^{i - j}N_{\alpha_{n-1}+}\lp -1,0\rf  \leq \gamma\lp j\rp, j= 0,\dots,i-1 \},\label{alphaad}
  \end{align}
  where $N$ is the $K$-poisson process driven by $\pi^{\downarrow \underline{\lambda},\overline{\lambda}}$. For setup (O), let $\varsigma_{\tau_{n}-\alpha_{n-1}}$ be the Dirac-measure on $\tau_{n}-\alpha_{n-1}$ and let $Z^{n,pre}$ be the linear Hawkes process driven by $\pi_{\alpha_{n-1}+}^{\downarrow \underline{\lambda},\overline{\lambda}}$ with weight function $h_{ + }$, rate function $\psi_{L}$ and initial signal \\ $R_{t}=f\lp t\rp+h_{ + }\lp t-(\tau_{n}-\alpha_{n-1})\rp\ql \lb t> \tau_{n}-\alpha_{n-1}\rb$. Define $Z^n=\varsigma_{\tau_{n}-\alpha_{n-1}}+Z^{n,pre}$ and put
{\small
  \begin{align*}  \alpha_{n} = \alpha_{n-1}+
\inf\left\{ i>\lceil \tau_{n}-\alpha_{n-1}\rceil: \int_{0}^{i}h_{ + }\lp t - s\rp  dZ^{n}_{s}  \leq f(t,i-1)\  \forall t>i \text{ and } Z^n\lp 0,i \rf\leq \int_{0}^i \gamma\lp s+1\rp ds \right\}.
  \end{align*}
}
\begin{rem} We notice some properties of the system above.
 \begin{enumerate}
\item $Z$ is a well-defined $\cF _{t}$-progressive process on $\R_{ + }$ and $\lambda,\underline{\lambda},\overline{\lambda}$ are $\cF_{t}$-predictable.

\item The system remains unchanged for any choice of initial conditions $R^*,A_{0}^*$ as long as $\alpha_{0}$ from \cref{As:alpha0} remains unchanged.
\item By \cref{rooflm} $\pi^{\downarrow \underline{\lambda},\overline{\lambda}},\pi^{\uparrow \underline{\lambda},\overline{\lambda}}$ are $\cF_{t}-PRMs$.
 In particular $\lp \alpha_{n}-\alpha_{n-1},\tau_{n}-\alpha_{n-1}\rp\indep \cF_{\alpha_{n-1}}$ given $\alpha_{n-1}<\infty.$
\item The process is reversible in the sense that  $\underline{\lambda}_{\alpha_{n}+s},\overline{\lambda}_{\alpha_{n}+s}$ may be computed from\\ $\lp \pi_{\alpha_{n}+}^{\downarrow \underline{\lambda},\overline{\lambda}}\rp_{\vert \lp 0,s \rf},\lp \pi_{\alpha_{n}+}^{\uparrow \underline{\lambda},\overline{\lambda}}\rp_{\vert \lp 0,s \rf}$. That is, there is a map $H:M^c_{\R_{+}^2}\times M^c_{\R_{+}^2} \times \R_{+}\con M^c_{\R_{+}^2}$ satisfying
\begin{align}\label{inversfunc}
H\lp \lp \pi_{\alpha_{n}+}^{\downarrow \underline{\lambda},\overline{\lambda}}\rp_{\vert \lp 0,s \rf},\lp \pi_{\alpha_{n}+}^{\uparrow \underline{\lambda},\overline{\lambda}}\rp_{\vert \lp 0,s \rf} ,s \rp=  \lp \pi_{\alpha_{n}+}\rp _{\vert \lp 0,s\rf},
\end{align}
for all $n\in \N_{0}$ s.t. $\alpha_{n}<\infty$.
\end{enumerate}
\end{rem}
The rest of this section is dedicated to presenting our main result \cref{regenerationTheorem}, and its related results \ref{Prop:AlphaExists} - \ref{s:AlphaDist}. Before we state it, we colloquially explain the essence of the result. The purpose of the $\alpha$'s is to have points in time, where the intensity contribution from the past of $Z^{*}$ may be replaced by something deterministic. More preciesly $\alpha_{n}$ should satisfy
	\begin{align}
  &\lv \int_{0}^{\alpha_{n}}h\lp t - s\rp  dZ^*_{s}+ R^{*}_{t} \rv  \leq f(t-\alpha_{n}),\quad t>\alpha_{n},n\in\N: \alpha_{n}<\infty.\label{alphan}
	\end{align}
We prove that this property holds for both setups in \cref{Prop:AlphaExists}. The inequality $\eqref{alphan}$ combined with the properties of $\psi$
gives that $\lambda\leq \lambda^*,$ at least locally in time after $\alpha_{n}$.  We will also be able to control the difference $\lambda^*-\lambda$ locally after $\alpha_n$, and establish that $Z$ mimics $Z^{*}$ in that same interval. In fact, the purpose of $\tau_{n+1}$ is to act as a conservative right end of an interval starting at $\alpha_n$ on which $Z$ and $Z^*$ are equal.  All this will be proved in \cref{s:bandlemma}.
We then proceed to study the distributions of $\tau_{n},\alpha_{n}.$ In  \cref{s:distributiontheorem} we study $\tau_{n+1}-\alpha_{n}\vert \alpha_{n}<\infty$ and prove $P(\tau_{n+1}=\infty,\alpha_{n}<\infty)>0$ along with moment properties of the distribution $\tau_{n+1}-\alpha_{n}\vert \tau_{n+1}<\infty$. In \cref{s:AlphaDist} we investigate the law of $\alpha_{n}-\tau_{n}\vert \tau_{n}<\infty$. Here we prove that $P\lp \alpha_{n}-\tau_{n}<\infty\vert \tau_{n}<\infty\rp =1, $ and we characterize its moments. Combining these results implies that
\begin{align}{\label{eq:eta}
	\eta := \inf\left\{ n \in \N_{0}: \tau_{n + 1} = \infty\right\}
}\end{align}
is finite almost surely, and we will be able to show that
\begin{align}
\rho:=\alpha_{\eta}+D
\end{align}
 is a point of regeneration for $Z^*$. In fact, we have $Z^*\lp \rho-D,\rho\rf=0$ so $Z^*_{\vert \lp 0,\rho\rf}\indep Z^*_{\alpha_{\eta}+}$, i.e. there is an overlap of length $D$. We characterize the integrability properties of $\rho$ by applying the previously mentioned results concerning $\alpha_{n},\tau_{n}$.
\begin{figure}[H]
	\centering
	\begin{tikzpicture}[scale=0.8]

	\draw[->] (-3.2,0) -- (10,0);
	\draw[->] (-3,-0.2) -- (-3,7);


	\draw (-3.5, 2.5) node {$\lambda^{*}$};
	\draw[rounded corners=4pt]
	(-3,2.5) to[bend left]
	(-2,3) to[bend left]
	(-1,2.2) to[bend left]
	(0,3) to[bend left]
	(0.5, 3.7) to[bend left]
	(1,3.2) to[bend left]
	(1.5,3.5) to[bend left] (2,4.5);
	\filldraw (2,4.5) circle (2pt);

	\draw (0.0,-0.3) node {$\alpha_0$};

	\draw[draw=blue] (0,0) circle (2pt);
	\draw (-0.4,0.4) node {\color{blue}$\lambda$};
	\draw[draw=blue] (0,0) to[bend left] (2,3);
	\filldraw[fill=blue] (2,3) circle (2pt);

	\draw (2,1) node[regular polygon,regular polygon sides=3, draw, scale=0.3, fill=black] {};

	\draw (-1, 6.6) node {\color{red}$\lambda+F(t)$};
	\draw[draw=red] (0,6.8) to[bend right] (2,6.3);
	\filldraw[fill=red] (2,6.3) circle (2pt);

	\draw[fill=red,nearly transparent] (0,6.8) to[bend right]
	(2,6.3) --
	(2,3) to[bend right]
	(0,0) -- cycle;


	\draw[draw=blue] (2,2) ..controls (2.5,1.5) and (3.2,3.8).. (4,3.8);
	\draw[draw=blue] (2,2) circle (2pt);
	\filldraw[fill=blue] (4,3.8) circle (2pt);

	\draw[rounded corners=4pt]	(2,3.7) to[bend left]
	(2.5, 4.2) to[bend left]
	(3,3.9) to[bend left]
	(3.5, 4.5) to[bend left]
	(4, 4.5);
	\draw (2,3.7) circle (2pt);
	\filldraw (4,4.5) circle (2pt);

	\draw[draw=red] (2,5.3) to[bend right] (4,5);
	\filldraw[fill=red] (4,5) circle (2pt);
	\draw[draw=red] (2,5.3) circle (2pt);

	\draw[fill=red, nearly transparent]
	(2,2) ..controls (2.5,1.5) and (3.2,3.8).. (4,3.8) --
	(4,5) to[bend left] (2,5.3) -- cycle;

	\draw (4,4.15) node[regular polygon,regular polygon sides=3, draw, scale=0.3, fill=black] {};

	\draw[-] (4,-0.2) -- (4,0.2);
	\draw (4,-0.4) node {$\tau_1$};


	\draw (4,2.3) circle (2pt);
	\draw[rounded corners=4pt]
	(4,2.3) to[bend left]
	(5,2.9) to[bend left]
	(6,2.5) to[bend left]
	(7,3) to[bend left]
	(7.5, 2.7) to[bend left]
	(8,3.5) to[bend left]
	(8.5,3.7) to[bend left] (9,4.5);
	\draw[dashed] (9,4.5) -- (10,4.9);

	\draw[draw=blue] (7,0) circle (2pt);
	\draw[draw=blue] (7,0) to[bend left] (9,3);
	\draw[dashed,blue] (9,3) -- (10,3.25);
	\draw (7,-0.4) node {$\alpha_1$};

	\draw[draw=red] (7,6.8) to[bend right] (9,6.1);
	\draw[dashed,draw=red] (9,6.1) to (10,6.2);

	\draw[fill=red,nearly transparent] (7,6.8) to[bend right]
	(9,6.1) --
	(9,3) to[bend right]
	(7,0) -- cycle;

	\end{tikzpicture}
	\caption{An illustration of the system with $D=0$. The points of the PRM $\pi$ is depicted by $\lp  \blacktriangle \rp$, and the three intensities $\underline{\lambda},\lambda^*,\overline{\lambda}$. is colored in blue, black and red respectively. The red band illustrates the area where $\pi$ and $\pi^{\downarrow \underline{\lambda},\overline{\lambda}}$ differs.}
	\label{fig:renewal}
\end{figure}
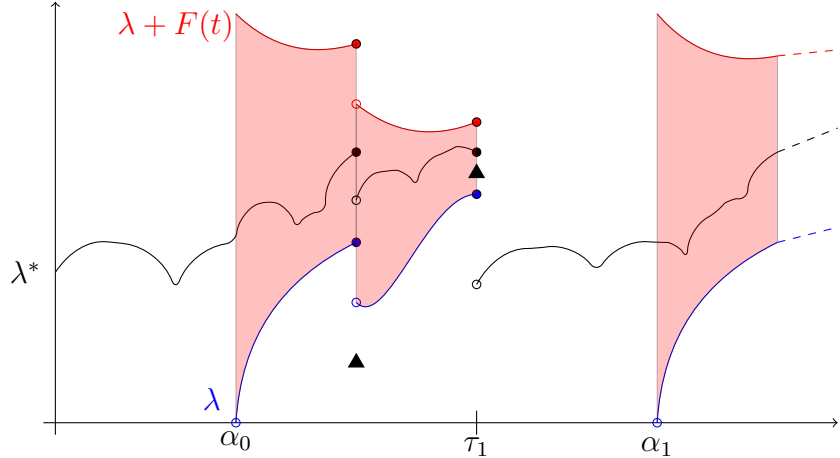
 The precise result is as follows:
Define $\cF^{*}_{t} = \sigma\lp \cF_{t}, \pi^{\uparrow \underline{\lambda},\overline{\lambda}} \rp $. It is clear that $(\cF^*_t)$ defines a filtration and without changing notation, we extend it to satisfy the usual hypothesis.

\begin{s}\label{regenerationTheorem}
	Grant either setup (AD) or (O).
	\begin{enumerate}
		\item  The random time $\rho = \alpha_{\eta}+D$ is an a.s. finite $\cF^*_t$-stopping time.
		\item The random measure $\pi^{\downarrow \underline{\lambda},\overline{\lambda}}$ is an $\cF^*_t$-PRM and hence $ \pi^{\downarrow \underline{\lambda},\overline{\lambda}}\indep \cF^*_{\rho}$. Moreover, we have $Z^*_{\rho+}=Z_{\rho+}$ and independent of $Z^{*}_{\vert \left(0,\rho\right] } $. In particular, $Z^*_{\rho+}$ is distributed as an ADHP with weight $h$, rate $\psi$, initial age $D$ and signal $t\mapsto - f\lp t+D\rp$.
		\item It holds that $\qE \lp \rho-\alpha_{0}\rp^p<\infty$. Under assumption 2B  it holds that $\rho-\alpha_{0}$ has exponential moment.
	\end{enumerate}
\end{s}
To prove \Cref{regenerationTheorem} we establish the results below and combine them in the end. The proofs of these results, and the main result, may be found in the proof section.

%

\begin{s}\label{Prop:AlphaExists}
 \textcolor{white}{h}\\[-35pt]
 \begin{enumerate}
 \item Consider setup (AD).
   It holds that $\alpha_{n}$, $n\in \N$ satisfies \eqref{alphan}. Moreover, assume $Z^I$ from \cref{eks:inv} exists, and assume only that \cref{As:fmap} (A) holds with $r=f$.  Then \cref{As:alpha0} is satisfied for $Z^*=Z^I$ and
\begin{align}
   \alpha^I_{0}&= \inf\lb i > 0 :\theta^{i-j}N\lp  -1 , 0  \rf \leq \gamma\lp j\rp , j\geq 0 \rb.
\end{align}
   \item Consider setup (O) and recall $\psi_{L}$ from \eqref{linear}.
   It holds that $\alpha_{n}$, $n\in \N$ satisfies \eqref{alphan}.  Moreover, assume $Z^I$ from \cref{eks:inv} exists, and assume only that \cref{As:fmap} (A) holds with $r=f$. Define $Z^0$ as the stationary linear Hawkes process driven by $\pi$ with weight/rate  $h_{+},\psi_{L}$ (see \cite{bm} theorem 1 and remark 8). Then  \cref{As:alpha0} is satisfied with
{\small
  \begin{align}
\label{alphao}   \alpha_{0}^I =
\inf\{ i>0: \int_{-\infty}^{i}h_{ + }\lp t - s\rp  dZ^{0}_{s}  \leq f(t-i)\;\; \forall t>i\}.
\end{align}
}
 \end{enumerate}
\end{s}

\begin{prop}\label{s:bandlemma}
Under either setup  it holds a.s. for all $n\in\N$ such that $\alpha_n<\infty$ and $t\in (\alpha_n, \tau_{n+1}]$ that
\begin{align}
Z(\alpha_{n}, \tau_{n+1}) = Z^*(\alpha_{n}, \tau_{n+1})
\end{align}
and
\begin{align}
0 \leq \lambda^*_t - \lambda_t \leq F(t-\alpha_n).
\end{align}
\end{prop}

\begin{s}\label{s:distributiontheorem}
 Under either setup it holds that
 \begin{align}\label{taudist}
   P\lp \tau_{n } - \alpha_{n - 1}=\infty \vert \alpha_{n-1}<\infty \rp &=\exp\lp -\lV F \rV_{\cL^{1}} \rp\\
P\left( \tau_{n } - \alpha_{n - 1}\leq t\vert  \tau_n<\infty  \right)  &=  \dfrac{1 - \exp\left( -\int_{0}^{t}F\left( s\right) ds \right) }{1-\exp\left(- \lV F \rV_{\cL^{1}} \right)}.\label{taudist2}
 \end{align}
 In particular the conditional distribution $\tau_{n } - \alpha_{n - 1}\vert \tau_{n } < \infty  $ has $p$'th moment. Under \cref{As:fmap} (B) it has exponential moment.
\end{s}

It turns out that the $\alpha_{n}$'s defined above may be analyzed using a discrete Markov chain. In fact one may rewrite $\alpha^I_{0}$ and $\alpha_{n}-\tau_n$ as return times to state 0 for a specific Random Exchange process (see Appendix). This yields precise distribution results as given in the next proposition.

\begin{s}\label{s:AlphaDist}
Under either setup it holds that $\alpha_{0}^I$ and $(\alpha_n - \tau_n)\mid (\tau_n<\infty)$ have $p$'th moment. If also \cref{As:fmap} (B) holds then these laws  have exponential moment.

%
%
%
\end{s}

\subsection{Hawkes Processes in a Markov Chain Framework}\label{markovsection}
In this section we first apply \cref{regenerationTheorem} iteratively to obtain consecutive renewal time points $\rho_{i}$, which partition $Z^*$ into independent bits. Afterwards, we construct a Markov chain that contains the information of $Z^*$, and where the $\rho_{i}$'s acts as the return times to an atom. The purpose is to use Markov chain theory to obtain results for $Z^*$, which we do in the next section.\\

Choose $D>0$ and set $\rho_{0}:=\rho,\pi^0:=\pi,\overline{\pi}^0:=\overline{\pi}$ and $\pi^1 =\pi^{\downarrow \underline{\lambda},\overline{\lambda}}_{\rho_{0}+}$. We also introduce another auxiliary PRM $\overline{\pi}^1$ independent of $ \pi^0,\overline{\pi}^0,\pi^1$. Note that $\pi^1,\overline{\pi}^1$ are $\cF^1_{t}$-PRMs where $\cF^1_t = \sigma(\cF^*_{\rho}, \pi^1_{\mid\lp 0,t\rf}, \overline{\pi}^1_{\mid\lp 0,t\rf})$. It holds that $Z^{*}_{\rho+}$ is an ADHP driven by $\pi^1$ with $R:t\mapsto -f\lp t+D\rp $, $A_{0}=D$. Thus with these PRMs $\alpha_{0}=0$ satisfies \cref{As:alpha0} with $r=f$. Thus, $\pi^1, \overline{\pi}^1$ and $\alpha^1_0$
induces a new renewal system as in section 3. From there we obtain sequences $\lp \alpha^{1}_{i}\rp,\lp \tau^{1}_{i}\rp$ and two new PRMs which we denote $\pi^{\downarrow 1},\pi^{\uparrow 1}$. Applying \cref{regenerationTheorem} on this system gives a new renewal time $\rho_{1}$.
\\\\
Continuing this way gives sequences $\lp\pi^i\rp,\lp \overline{\pi}^i\rp$,$\lp\pi^{\uparrow i}\rp,\lp\pi^{\downarrow i}\rp$, $\lp \rho_{i}\rp $. If we set $\varrho_{i}=\sum_{j=0}^i \rho_{j}$ we have that
\begin{align}
B\mapsto Z^{*}\lp \lp \varrho_{i-1}-D,\varrho_{i}\rf\cap \lp B+\varrho_{i-1}\rp  \rp, \quad B\in \cB_{\lp -D,\infty\rp }.
\end{align}
are i.i.d for $i\in\N$, each being the Hawkes process initialized with $A_{0}=D,R:t\mapsto -f\lp t+D\rp $ and driven by $\pi^i.$  \\\\ In fact, we can study $Z^{*}$ from a Markov chain perspective.
Define
\begin{align}
\rho^{-}\lp t\rp&=\sup \lb s<t:\exists i\in \N_{0} : s=\varrho_{i}\rb,\\
J\lp n\rp &=\lv \lb j\in \N_{0}:\varrho_{j}< n  \rb\rv,
\end{align}
and $A^{\rho}_{t}=t-\rho^{-}\lp t\rp $. Consider the stochastic processes on the state-space
$M^{c}_{\R^2_{+}}\times M^{c}_{\R^2_{+}}\times \N$,
\begin{align} \label{eq:markovlol}
\Phi^{pre}_{n}&=\lp \pi^{\downarrow J\lp n\rp}_{\vert \lp 0,A^{\rho}_{n}\rf },\pi^{\uparrow J\lp n\rp },A^{p}_{n} \rp\\
\Phi_{n}&=\theta^{\rho_{0}}\Phi^{pre}_{n}
\end{align}
for $n\in \N$. In this framework, \Cref{regenerationTheorem} states that $ \Phi_{n}  \indep \Phi^{pre}_{0},\dots , \Phi^{pre}_{\rho_{0}}$. Using \eqref{inversfunc}  we may construct a map $H^*:M^{c}_{\R^2_{+}}\times M^{c}_{\R^2_{+}} \times \N\con M^{c}_{\lp -D,0\rf}$ satisfying\\ $H^*\lp \Phi_{n}\rp = \lp \theta^{\rho_{0}+n }Z^*\rp_{\vert \lp -D,0 \rf }$.
  Also, by construction of $\rho_{i}$, the indicator function \\$J\lp n+\rho_{0}+1\rp- J\lp n+\rho_{0}\rp$ may be written as some map $H_{J}\lp \Phi_{n}\rp$.  It follows that $\Phi$ is a Markov chain with an atom
\begin{align}
\Xi = \lb \mu,\nu,n: H_{J}\lp \mu,\nu,n\rp=1\rb.
\end{align}
Consider the subspace
\begin{align}
\qX = \lb \mu,\nu,n: \mu,\nu \text{ are simple }, P^\Phi_{\mu,\nu,n}\lp \Phi \text{ hits } \Xi \text{ eventually}\rp =1\rb
\end{align}
and let $P^{\Phi}$ be the kernel of $\Phi$. By definition of $\qX$, any chain with kernel $P^{\Phi}$ started in $\qX$ eventually hits $\Xi$, and by \cref{regenerationTheorem} $\Phi$ almost surely returns to $\Xi$ once hitting it. Thus $\qX$ is an absorbing state for chain $\Phi$, and we shall from now on always refer to $\Phi$, $P^{\Phi}$ as the restricted kernel to $\qX$ (see proposition 4.2.4. \cite{MEYN}). Since $\Xi$ is an accesible atom for this chain, it is irreducible (Prop. 5.1.1 \cite{MEYN}) and aperiodic. The return time to $\Xi$ is distributed as $\rho_{1}$, so by Kac's theorem it follows that $P^{\Phi}$ is positive with invariant law $\tilde{P}$ for $p\geq 1$. Not surprisingly it turns out that $\tilde{P}$ agrees with $Z^I$ from \cref{eks:inv} with parameters $h,\psi$, whenever they both exist.

\begin{prop}
\label{prop:inv} Assume that $p\geq 1$ and $Z^I$ from  \cref{eks:inv} exists with coupling properties as given in the example. It holds that $H^*\lp \tilde{P} \rp \stackrel{D}{=}Z^I_{\vert \lp -D,0\rf}$.
\end{prop}

\section{Applications}
In this section we apply the Markov chain construction from section \ref{markovsection} to establish asymptotic results for $Z^*.$ We show distribution results of the coupling time. Then we present a functional CLT, a time-average CLT, and a LIL for $Z^*$.

\subsection{Bound On The Coupling Time}

In \cite{bm} it was shown that two ordinary Hawkes processes, started with different initial conditions couple under regularity conditions. In \cite{MES} we showed a similar statement for ADHPs. our construction of $\rho$ confirms these results, and provide moment results for these coupling times. More preciesly, we have the following proposition
\begin{prop}\label{CouplingBound}
	Let $Z^{*1},Z^{*2}$ be two $\pi$-driven, ADHPs initialized with $A^{*1}_0,A^{*2}_0$ and signals $R^{*1}_t, R^{*2}_t$. Assume that $\alpha_{0}$ satisifies \cref{As:alpha0} for both measures simultaneously.\\ Define the coupling time  $T=\inf\lb t>0: \lv Z^{*1}-Z^{*2}\rv \lp t,\infty\rp = 0\rb $.
It holds that $T\leq \rho$ and in particular that $\lp T-\alpha_{0}\rp_{+}$ has p'th moment. Under \cref{As:fmap} $\lp B\rp$, $\lp T-\alpha_{0}\rp_{+}$ has exponential moment.
\end{prop}
\begin{p}
By \cref{regenerationTheorem} $2)$ we have $Z_{\rho+}^{*1}=Z_{\rho+}^{*2}=Z_{\rho+}$. The moment results follows from \cref{s:distributiontheorem} and \cref{s:AlphaDist}.
\end{p}
\subsection{Asymptotics}
%
%

The Markov chain $\Phi_{n}$ from \eqref{eq:markovlol} can be used to establish various asymptotic results for $Z^*$ in a general setting. Let $G:M^c_{\R_+^2}\times M^c_{\R_+^2} \times \N \con \R$ be a measurable function which we normalize with  $\overline{G}=G-\tilde{P}G$ where $\tilde{P}$ is the invariant measure from \cref{prop:inv}. We shall discuss asymptotic results of the sum
\begin{align}
S_{n}\lp \overline{G}\rp=\sum_{k=1}^n \overline{G}\lp \Phi^{pre}_{k}\rp
\end{align}
provided of course that $G$ is a function s.t. $G\lp \Phi^{pre}_{n}\rp $ is a well-defined variable for all $n\in \N$. Define
\begin{align}
\tilde{S}_{n}\lp \overline{G}\rp=S_{n+\rho_{0}}-S_{\rho_{0}}.
\end{align}
\begin{s}
\label{th:CLTaverage}
 Assume that $p\geq 2$. Define $\mu_{\rho} = \qE \rho_{1}$ and
\begin{align}
\sigma^2=\mu_{\rho}^{-1} \qE \tilde{S}_{\rho_{1}}\lp \overline{G}\rp^2.
\end{align}
Assume that $\sigma^2$ is finite and nonzero.
\begin{enumerate}
\item The following CLT holds
\begin{align}
n^{-1/2}S_{n}\lp \overline{G}\rp  \Rightarrow N\lp 0,\sigma^{2}\rp.
\end{align}

\item The following LIL holds: Almost surely,

\begin{align}
\liminf_{n\con \infty}\frac{S_{n}\lp \overline{G}\rp }{ \sqrt{2\sigma^2 n\ln\ln n}}=-1,\quad \limsup_{n\con \infty}\frac{S_{n}\lp \overline{G}\rp }{ \sqrt{2\sigma^2 n\ln\ln n}}=1.
\end{align}
\end{enumerate}

\item Define $S_{t}\lp \overline{G}\rp$ for $t\in \lp n,n+1\rp$ as the linear interpolation between $S_{n}\lp \overline{G}\rp$ and $S_{n+1}\lp \overline{G}\rp$, and put
\begin{align}
B^n_{t}= \frac{1}{\sqrt{n\sigma^2}}S_{nt}\lp \overline{G}\rp.
\end{align}
The following functional CLT holds
\begin{align}
B^n_{\lp \cdot \rp }\Rightarrow  B_{(\cdot)}
\end{align}
where $B$ is the 1-dim Brownian motion, and the convergence ($\Rightarrow$) is in distribution in the D space on $\lf 0,\infty \rp$. See \cite{billingsley} chapter 3.
\end{s}
\begin{p}
A direct application of theorem 17.2.2, theorem 17.4.4, and section 17.4.3 \cite{MEYN}  on $\Phi_{n}$  gives the desired results for $\tilde{S}_{n}$  inplace of $S_{n}.$ Standard arguments extends these results to hold for $S_{n}$ as well.
\end{p}

\begin{eks} \label{eks:intex}
We give two applications of \cref{th:CLTaverage}. \\\\
\textbf{1. } Take $G\lp \mu\rp=\int_{-D}^0 \phi\lp s\rp  dH^*\lp \mu\rp\lp s\rp  $ where $\phi:\lp -D,0\rf \con \R$ is a bounded map. In particular, when $\phi\equiv 1,D=1,$ we obtain
 \begin{align}
S_{n}\lp \overline{G}\rp= Z^*\lp 0,n\rf-\qE Z^I\lp 0,n\rf.
\end{align}
 It is straightforward to show a sufficient criteria for $\sigma^2<\infty$ is that either $p > 2$ in the (O) setup or $p\geq 2$ in the (AD) setup. Also it is easy to see that for non-degenerate choices of $h,\psi$ we have $\sigma^2>0$.
\\\\
\textbf{2. }
 Fix $m\in\N$, and let $T:M^c_{\vert \lp -m,0\rf}\con \R$ be a measurable map. Take $D=m+1$ and $G\lp \mu\rp =\int_{0}^1 T\lp \lp \theta^s  H^*\lp \mu \rp\rp_{\vert \lp -m,0\rf}\rp  ds $. These choices leads to
\begin{align}
S_{n}\lp \overline{G}\rp=\int_{0}^{n} T\lp  \theta^{s} Z^* _{\vert \lp s-m,s\rf}\rp ds -  n\qE T\lp Z^I_{\vert \lp -m,0\rf}\rp.
\end{align}

To obtain $\sigma^2<\infty$ we need a growth condition  depending on the setup:

\begin{align}
\text{(AD)}&\quad  \exists c,C>0: \lv T\lp \mu \rp \rv\leq C\exp \lp c\mu\lp -m,0\rf\rp \\
\text{(O)}&\quad \exists C>0\quad:    \lv T\lp \mu \rp \rv \leq C\lv \dfrac{\mu\lp -m,0\rf}{1+\ln\lp \mu\lp -m,0\rf\rp }\rv^{p/2-1 }.
\end{align}
Consider first setup $(AD)$. Recall that
\begin{align}
  Z_{\rho_{0}+}^*\lp n-D,n\rf=Z_{\rho_{0}+}\lp n-D,n\rf\leq D\delta^{-1}+\pi^{\downarrow \underline{\lambda},\overline{\lambda}}_{\rho_{0}+}\lp n-D,n\rf\times \lf 0,K\rf.
\end{align}
 It follows that
\begin{align}
\sup_{s\in\lp n-1,n\rf} \lv T\lp \theta^{-s} Z^*_{\vert \lp s-m,s \rf}\rp \rv \leq C\exp\lp cD\delta^{-1}\rp  \exp\lp c \pi^{\downarrow \underline{\lambda},\overline{\lambda}}_{\rho_{0}+}\lp n-D,n\rf\times \lf 0,K\rf \rp :=Y_{n}
\end{align}
For $j=1\dots D,n\in \N$ define $\tilde{Y}^j_{n}=Y_{n}$ if $n\equiv j$ mod $D$, and otherwise set $\tilde{Y}^j_{n}$ as an i.i.d copy of $Y_{1}$. Indeed by the $C_{p}$-inequality there is some possibly larger $C>0$ s.t.
\begin{align}
\tilde{S}_{\rho_{1}}\lp G\rp^2\leq C\sum_{j=1}^D \lp \sum_{n=1}^{\rho_{1}} \tilde{Y}^j_{n}\rp^2
\end{align}
Notice that for each fixed $j=1,\dots,D$, the sequence $\tilde{Y}^j_{n}, n\in \N$ are i.i.d. It follows from  theorem 5.2, chapter 1 in \cite{Gut} that $\sigma^2<\infty$.
\\\\
Consider now setup $(O)$, and take $\gamma\lp t\rp =C\ln_{+}t $ for $C$ so large that \eqref{criteriaO} is satisfied. Note that $x\mapsto x/\ln\lp x+1\rp $ is increasing for $x>0$ so
\begin{align}
\sup_{s\in\lp \varrho_{0},\varrho_{1}\rf} \lv T\lp \theta^{-s} Z^*_{\vert \lp s-m,s \rf}\rp \rv &\leq C\lv \dfrac{Z^*\lp \varrho_{0},\varrho_{1}\rf}{\ln\lp 1+Z^*\lp \varrho_{0},\varrho_{1}\rf \rp }\rv^{p/2-1 }.
\end{align}
Notice that for our choice of $\gamma$, we have the inequality $\int_{0}^x\gamma\lp t+1 \rp dt\leq C x\ln \lp x+1\rp$ for a possibly larger constant $C>0$. From the definition of $\alpha_{n}$ it follows that

\begin{align}
\sup_{s\in\lp \varrho_{0},\varrho_{1}\rf} \lv T\lp \theta^{-s} Z^*_{\vert \lp s-m,s \rf}\rp \rv\leq  C_{T} \lv \dfrac{C \rho_{1}\ln \lp \rho_{1}+1\rp}{\ln \lp 1+ C \rho_{1}\ln \lp \rho_{1}+1\rp\rp   }\rv^{p/2-1 }\leq C_{T}C \rho_{1}^{p/2-1}
\end{align}
again for a possibly larger constant $C>0$. From here it follows that $\sigma^2<\infty$. One would have to check $\sigma^2>0$ for the given $T$, but for most practical applications, this is a triviality.

\end{eks}

\section{Discussion and Outlook}
In the following, we shall discuss generalizations and limitations of the presented results, and suggest further research topics.
\subsubsection*{Multivariate Hawkes processes}
It is straight forward to generalize the regeneration procedure to a multivariate Hawkes Process with $N$ units (see \cite{dfh},\cite{SusEva} or \cite{MES} for an introduction to these). One should split each $\pi_{i}$ , $i\leq N$  into $\pi^{\uparrow i},\pi^{\downarrow i}$
for $i=1,...,N$ - analogous to what was done in the start of section 3.1. The $\tau_{n}^i$'s should be generalized in the obvious way, while $\alpha_{n}$ should be modified so that it ensures that $\sum_{i=1}^N
\lv \int_{0}^{\alpha_{n}} h_{ij}\lp t-s\rp dZ^i_{s}+ R^i_{t}\rv
\leq f\lp t-\alpha_{n}\rp.$
In setup (AD) this is achieved by substituting $\pi^{\downarrow \underline{\lambda},\overline{\lambda}}$  in \eqref{alphaad},\eqref{alphao} with $\sum_{i=1}^N \pi^{\downarrow i}$ which will be a PRM with mean intensity $ N \; dzds$. in setup (O), the clusters $Z^{i}$ should be dominating linear $N$-dimensional Hawkes processes. While the total progeny distribution of $Z^i$ is no longer Borel distributed, it is well known that it has exponential moment, which is sufficient to complete the proof.

\subsubsection*{Stability for more general setups}
A significant observation is that the setups (AD) and (O), essentially only affects the construction of $\rho$ through the choice of $\alpha$'s and $f$. For other and more general setups in the univariate or multivariate case, one may adapt this procedure to establish stability regimes. For example it might be a method to explore other multivariate systems where inhibition from either the weight or the age have a potential effect on the stability regime.

\subsubsection*{Optimizing the regeneration scheme}
These results establish a regeneration scheme for weight functions $h$ s.t.\\ $\int_{0}^\infty \lv h\lp t\rp\rv t^{p+1}\ln_{+}t dt<\infty$. However,  invariant solutions to $Z$ exists already for $h$ with first moment i.e. $\int_{0}^\infty t \lv h\lp t\rp \rv dt<\infty.$ Also, the CLT result for ordinary Hawkes processes by Zhu \cite{Zhu1}, assuming that $h$ decreasing and positive, only requires that $t\mapsto th\lp t\rp$ is integrable. This corresponds to $p=0$, instead of $p\geq 2$ which we require in \cref{th:CLTaverage}. These facts indicate that there may exist renewal times with better moment properties, than those discussed in this article.

\subsubsection*{Implementation and practical computation}
While this article focus on the theoretical development of regeneration times, the method is constructive and $\rho$ may be simulated in either setup. It would be of interest to study the efficiency of this algorithm.


\section{Proofs}

\subsection{Proofs of Section 2 Results}
\noindent\textbf{Proof of \cref{prop:fmap}}\\
To show the claimed result for $f$ under \cref{As:fmap} (A), substitute the inner variable with $u=s+t$ and apply Tonelli to obtain
\begin{align}
\int_{0}^{\infty} \int_{0}^{\infty} t^{p}\gamma\lp s+1\rp  \overline{h}\left( t + s\right) ds dt &= \int_{0}^{\infty} \int_{t}^{\infty} t^{p}\gamma\lp u-t+1\rp  \overline{h}\left( u\right) dudt\\
&= \int_{0}^{\infty} \int_{0}^{u} t^{p}\gamma\lp u-t+1\rp  \overline{h}\left( u\right) dtdu
 \\
&\leq\lp p+1\rp^{-1}\int_{0}^{\infty} u^{p+1}\gamma\lp u+1\rp \overline{h}\left( u\right) du
\end{align}
which proves the desired. Under \cref{As:fmap} (B) it is straightforward to show that  $f$ has exponential moment. The claimed result for $F^p,F$ follows immediately.
\subsection{proofs of section 3.1 results}
\noindent\textbf{Proof of \cref{Prop:AlphaExists}}\\

\textbf{1.}
 We claim for $n\in \N$ that
 	\begin{align}
 	\lv \int_{0}^{\alpha_{n}}h\lp t-s\rp dZ^{*}_{s}+R^{*}_{t}\rv  \leq f\lp t-\alpha_{n}\rp
\quad  \forall t>\alpha_{n} \label{bandInequalityAD}
 	\end{align}
 	if $\alpha_{n}<\infty$.
We prove this by induction over $n\in \N_{0}$. The induction start $n=0$ is per \cref{As:alpha0}. To prove the induction step, we split the integral of interest
 	\begin{align}
 	&\lv \int_{0}^{\alpha_{n}}h\lp t-s\rp dZ^{*}_{s}+R^{*}_{t}\rv  \leq \lv \int_{\alpha_{n-1}}^{\alpha_{n}}h\lp t-s\rp dZ^{*}_{s}\rv +\lv \int_{0}^{\alpha_{n-1}}h\lp t-s\rp dZ^{*}_{s}+R^{*}_{t}\rv.\label{eq:split}
 	\end{align}
  By the induction assumption, the 2nd term is bounded by $f\lp t-\alpha_{n-1}\rp$ for all $t>\alpha_{n - 1}$.
The first term above can be split up to whether jumps of $Z^*$ happen when $A^*_{t}\leq \delta$ or not
\begin{align}
\int_{\alpha_{n-1}}^{\alpha_{n}}\ql \lb A_{s} \leq \delta\rb h\lp t-s\rp dZ^{*}_{s} +\int_{\alpha_{n-1}}^{\alpha_{n}}\ql \lb A_{s} >\delta\rb h\lp t-s\rp dZ^{*}_{s}.\label{deltapart}
\end{align}
 	 Consider the first term. By the (AD) criteria we have
\begin{align}
   \int_{a}^{b} \ql \lb A_{s} \leq \delta\rb dZ^*_{s}\leq \pi^K\lp a,b\rf:=\int_{\lp a,b\rb \times \lf 0,K\rf } d\pi\lp s,z\rp.
\end{align}
 On the interval $t\in\lp \alpha_{n-1},\tau_{n}\rp $ we have per definition of $\tau_{n}$ that $N\lf t,t\rf\geq \pi^{K}\lf t,t\rf$ while for $t\in\lp \tau_{n},\alpha_{n}\rf$ we have $\pi^{K}\lf t,t\rf=N\lf t,t\rf$. It follows that for $t>\alpha_{n}$, we have

 	\begin{align}
 	&\int_{\alpha_{n-1}}^{\alpha_{n}}\ql \lb A_{s} \leq \delta\rb h\lp t-s\rp dZ^{*}_{s}\\
 \leq &\int_{\alpha_{n-1}}^{\alpha_{n}}\overline{h}\lp t-s\rp d\pi^K_{s}\\
 \leq &\overline{h}\lp t - \tau_{n}\rp  + \int_{\alpha_{n-1}}^{\alpha_{n}}\overline{h}\lp t-s\rp dN_{s}\\
 	\leq &\overline{h}\lp t - \alpha_{n}\rp  + \sum_{i = 0}^{ \alpha_{n} - \alpha_{n - 1}-1}\gamma\lp i\rp \overline{h}\lp t - \alpha_{n} + i\rp \\
\leq&\int_{0}^{\alpha_{n} - \alpha_{n - 1} - 1}\gamma\lp s + 1\rp \overline{h}\lp t-\alpha_{n} + s\rp   \;ds + \lp 1 + \gamma\lp 0\rp\rp \overline{h}\lp t-\alpha_{n}\rp.\label{alphabound}
 	\end{align}
For the second integral in \eqref{deltapart},  recall that $\delta^{-1}\in \N$ and notice that the interdistance between jumps of $Z^*$ with $A_{s}>\delta$ is at least $\delta$ per definition. We obtain the bound
 	\begin{align}
 &\int_{\alpha_{n-1}}^{\alpha_{n}}\ql \lb A_{s} >\delta\rb h\lp t-s\rp dZ^{*}_{s}\\
\leq &\sum_{i = 0}^{\lp \alpha_{n} - \alpha_{n - 1}\rp\delta^{-1}-1}\overline{h}\lp t - \alpha_{n} + i\delta\rp\\
\leq &\sum_{j = 0}^{\lp \alpha_{n} - \alpha_{n - 1}\rp-1}\delta^{-1}\overline{h}\lp t - \alpha_{n} + j\rp\\
\leq &\int_{0}^{\alpha_{n} - \alpha_{n - 1} - 1} \delta^{ - 1}\overline{h}\lp t-\alpha_{n} + s\rp   \;ds+ \delta^{ - 1}\overline{h}\lp t-\alpha_{n}\rp.\label{deltabound}
 	\end{align}
  The sum of the two right-hand sides of \eqref{alphabound},\eqref{deltabound}  are less than $f\lp t-\alpha_{n-1},\alpha_{n}-\alpha_{n-1}-1\rp.$ The induction claim now follows by inserting this back into \eqref{eq:split}.\\\\
 To prove that $\alpha_{0}^{I}$ satisfies \cref{As:alpha0}, repeat the proof above with $-\infty$ in place of $\alpha_{n-1}$ and $\alpha^{I}_{0}$ in place of $\alpha_{n}$. We omit the details.\\\\

\noindent\textbf{2.} Consider now the (O) setup. We claim again for $n\in \N$ that
 	\begin{align}
 	\lv \int_{0}^{\alpha_{n}}h\lp t-s\rp dZ^{*}_{s}+R^{*}_{t}\rv  \leq f\lp t-\alpha_{n}\rp
\quad  \forall t>\alpha_{n} \label{bandInequalityO}.
 	\end{align}
As before we proceed by induction over $n\in \N_{0}$. The induction start follows from \cref{As:alpha0}. Assume now the claim holds for $n-1$. By the induction assumption and the definition of $f$ we have

 	\begin{align}
 	\lv \int_{0}^{\alpha_{n}}h\lp t-s\rp dZ^{*}_{s}+R^{*}_{t}\rv  &\leq \lv \int_{\alpha_{n-1}}^{\alpha_{n}}h\lp t-s\rp dZ^{*}_{s}\rv + f\lp t-\alpha_{n-1}\rp. \label{thiso}
 	\end{align}
It is seen per induction over jumps of $Z_{n}$, that $Z_{n}\lf s,s\rf\geq Z^*\lf \alpha_{n-1}+s,\alpha_{n-1}+s\rf$ for all $s\in \lp \alpha_{n-1},\alpha_{n}\rf$. Hence

\begin{align}
\lv \int_{\alpha_{n-1}}^{\alpha_{n}}h\lp t-s\rp dZ^{*}_{s}\rv\leq \lv \int_{0}^{\alpha_{n}-\alpha_{n-1}}h\lp t-s\rp dZ^n_{s}\rv\leq f(t-\alpha_{n},\alpha_{n}-\alpha_{n-1}-1)
\end{align}
inserting back into \eqref{thiso}, and using the definition of $f$ gives the desired result. \\\\ The statement about  $\alpha_{0}^I$ is a direct implication of the claim that $Z^I\lf t,t\rf\leq Z^0\lf t,t\rf$ for all $t\in \R$ almost surely. To prove this claim, let $Z^1,Z^2$ be Hawkes processes with weight $h_{+},h$ and rate functions $\psi_{L},\psi$ respectively, with common intialization $R\equiv 0,A_{0}=0.$ By the coupling property from \cref{eks:inv} We have almost surely that $Z^0_{\vert \lf t,\infty\rp} =Z^1_{\vert \lf t,\infty\rp},Z^I_{\vert \lf t,\infty\rp} =Z^2_{\vert \lf t,\infty\rp}$ for $t$ large enough. On the other hand, per induction over jumps of $Z^1$ it is straightforward to prove that $Z^1\lf t,t\rf \geq Z^2\lf t,t\rf$ for all $t\in \R_{+}$ so it follows that almost surely $Z^0\lf t,t\rf \geq Z^I\lf t,t\rf $ eventually. The claim now follows from the fact that $\lp Z^0,Z^I\rp$ is stationary and ergodic.

\qed

\noindent\textbf{Proof of \cref{s:bandlemma}}\\

\noindent Define  $\tilde{Z} = \lv d\left( Z^{*} - Z\right) \rv$ and let $\tilde{A}$ be its age process. We have for $t\in \lp \alpha_{n-1},\tau_{n}\rf$
	\begin{align}\label{bandInequality2a}
	X^*_{t} - X_{t} & =  \int_{\alpha_{n - 1}}^{t - }  h\left( t - s\right) d(Z^*_s - Z_s)+\int_{0}^{\alpha_{n-1}}  h\left( t - s\right) dZ^{*}_{s}+R^{*}\lp t\rp  - f\lp t - \alpha_{n-1}\rp.
	\end{align}
	Applying \eqref{alphan} gives
	\begin{align}
	\lv \int_{0}^{\alpha_{n-1}}  h\left( t - s\right) dZ^{*}_{s}+R^{*}\lp t\rp  - f\lp t - \alpha_{n - 1}\rp\rv \leq 2f\lp t - \alpha_{n - 1}\rp \label{bandInequality2}.
	\end{align}
	Define for $n \in \N$
	\begin{align}
		\tau^{*}_{n} &= \inf\left\{ t>\alpha_{n - 1}:\tilde{Z}\left[ t,t \right] = 1\right\}.
	\end{align}
 We claim that $\tau^*_n \geq  \tau_n $ almost surely. Notice that $A^*_t \geq  A_t$ for $t\in (\alpha_{n-1}, \tau_n^*]$. For all $t\in\lp \alpha_{n-1},\lp \alpha_{n-1}+D \rp \wedge \tau^{*}_{n}\rf $ we combine \eqref{alphan} and \eqref{eq:psias2} to see that
	\begin{align}
	0 = \lambda_t \leq \lambda^*_t \leq c_{\psi}+Lf(t-\alpha_{n-1})\leq F(t-\alpha_{n-1})
	\end{align}
	which shows that $\tau_{n}^*\geq\tau_{n}\wedge \lp \alpha_{n - 1}+D\rp .$  For all $t$ in the (possibly empty) interval $\left( \alpha_{n-1}+D,\tau^{*}_{n} \right]$ we have
  \begin{align}
	X^{*}_{t} &=\int_{\alpha_{n-1}}^{t-}h\lp t-s\rp dZ^*_s+ \int_{0}^{\alpha_{n-1}}h\lp t-s\rp dZ^{*}_{s} + R^{*}_{t} \label{bandInequality2A}\\
&\geq \int_{\alpha_{n-1}}^{t-}h\lp t-s\rp dZ^{*}_{s}-f\lp t-\alpha_{n-1}\rp\\&=X_{t}.
	\end{align}
	Since $\psi$ is increasing we get the following inequality
	\begin{align}{
		\psi\left( X_{t},A_{t}\right) \leq \psi\left( X^*_{t},A^*_{t}\right) \quad \text{for all } t\in (\alpha_{n-1}+D, \tau_n^*].
	}\end{align}
	If $A^*_t>A_t$ the definition of $\tau_n^*$ gives that $A^*_t, A_t\geq t-\alpha_{n-1}$. Therefore  \eqref{bandInequality2a} and \eqref{bandInequality2} gives
	\begin{align}{
		\psi\left( X^*_{t},A^*_{t}\right) - \psi\left( X_{t},A_{t}\right) &\leq L\lp X^*_{t} - X_{t}\rp  +  c_{\psi}  g\left( t - \alpha_{n - 1}\right) \\
		& \leq  2L f\left( t - \alpha_{n - 1}\right)  +  c_{\psi} g\left( t - \alpha_{n - 1}\right)\\
		&\leq F\lp t-\alpha_{n-1}\rp.
	}\end{align}
	Likewise, if $A^*_t = A_t$ we have
	\begin{align}
	\psi\left( X^*_{t},A^*_{t}\right) - \psi\left( X_{t},A_{t}\right) \leq L  \lp X^*_{t} - X_{t}\rp
	\end{align}
	By definition of $\tau_{n}$, this implies  $\tau^{ *}_{n} \geq \tau_{n}.$
	Thus, in between two consecutive $\tau $ stopping times, the two Hawkes processes agrees.\qed

\vspace{1cm}

\noindent\textbf{Proof of \cref{s:distributiontheorem}}\\
The results \eqref{taudist},\eqref{taudist2} is a straightforward consequence of the strong Markov Property and \cite{bm} lemma 1. The density of the conditional distribution $\tau_{n}-\alpha_{n-1}\vert \tau_{n}<\infty$ is proportional to
\begin{align}
F\lp t\rp\exp\lp -\int_{0}^{t} F\lp s\rp ds \rp\leq F\lp t\rp
\end{align}
which shows the desired moment results for the distribution .\qed

\vspace{1cm}

\noindent\textbf{Proof of \cref{s:AlphaDist}}\\
To structure the proof, we discuss the following four variables, in written order
\begin{enumerate}
\item  $\alpha_{0}^I$ in the (AD) setup
\item $\alpha_{n}$ in the (AD) setup
\item $\alpha_{0}^I$ in the (O) setup
\item $\alpha_{n}$ in the (O) setup
\end{enumerate}
\newpage
  \textbf{1.}  To study $\alpha^I_{0}$, we introduce the process $M_{i}$ for $i\geq 0 $ given by
		\begin{align}
		M_{i} &= \inf\lb m \geq 0:\theta^{i-j}N\lp -1 ,0\rf \leq \gamma\lp j+m\rp ,j \geq 0 \rb \label{Minv}\\
&=\inf\lb m \geq 0:\theta^{k}N\lp -1 ,0\rf \leq \gamma\lp i+m-k\rp , k\leq i \rb.
		\end{align}

\begin{figure}[H]
\scalebox{.4}{
	\begin{tikzpicture}[line width=1mm,font={\fontsize{30pt}{12}\selectfont}]
	\draw[line width = 0.5mm] (3,0) -- (36,0);
	\filldraw (2.5, 0) circle (2pt);
	\filldraw (2, 0) circle (2pt);
	\filldraw (1.5, 0) circle (2pt);

	\draw (20,-0.5) node[below] {$i$} -- (20, 0.5);
	\draw (16,-0.5) node[below] {$i-1$} -- (16, 0.5);
	\draw (12,-0.5) node[below] {$i-2$} -- (12, 0.5);
	\draw (8,-0.5) node[below] {$i-3$} -- (8, 0.5);
	\draw (4,-0.5) node[below] {$i-4$} -- (4, 0.5);

	\draw (34, -0.5) -- (34, 0.5);

	\draw [blue, domain=2.5:30] plot (\x, {1.8 * ln(20 + 11 - \x)});
	\draw [blue] (30, 0) -- (34,0);

	\draw[ku] (18, 0)  -- (18, 2.5);
	\filldraw[ku] (18,0) circle (5pt);
	\filldraw[ku] (18,2.5) circle (5pt);

	\draw[ku] (14, 0)  -- (14, 4.8);
	\filldraw[ku] (14,0) circle (5pt);
	\filldraw[ku] (14,4.8) circle (5pt);

	\draw[ku] (10, 0)  -- (10, 1.5);
	\filldraw[ku] (10,0) circle (5pt);
	\filldraw[ku] (10,1.5) circle (5pt);

	\draw[ku] (6, 0)  -- (6, 4);
	\filldraw[ku] (6,0) circle (5pt);
	\filldraw[ku] (6,4) circle (5pt);

	\draw [ decorate,decoration={brace,amplitude=30pt,mirror},xshift=0.4pt,yshift=-0.4pt](20,0) -- (34,0) node[black,midway,yshift=-2cm] {$M_i$};

	\draw[blue] (27.3,5) -- (28.3,5) node[right] { $j\mapsto \gamma\lp i+M_{i}-j\rp $};

	\draw[ku] (30.3,3) -- (31.3, 3) node[right] { $N(i-1,i]$};
	\filldraw[ku] (31.3,3) circle (5pt);

	\end{tikzpicture}
}
\caption{ An illustration of $M_{i}.$ It describes how many right-shifts we need to apply to the curve $j\mapsto \gamma\lp i-j\rp $ before it bounds the entire partition $\lp N\lp i-j-1,i-j\rf\rp_{j\in \N_{0}}$}
\end{figure}
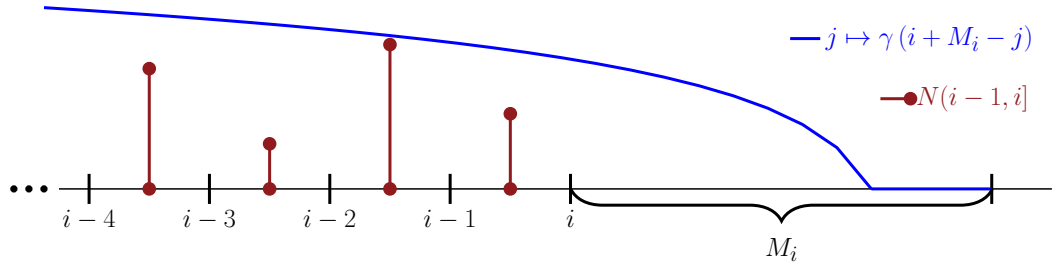

		We see that $\alpha_{0}=\inf\lb i>0:M_{i} = 0\rb$ i.e $\alpha_{0}$ is the first time $M_{i}$ hits $0$. Notice that $M_{0}$ is a well-defined random variable by a Borel-Cantelli argument. Observe also that $M$ is in fact an RE-process with update scheme
		\begin{align}
		M_{i} = \lp M_{i - 1}-1\rp \vee \lc \gamma^{ - 1}\lp N\lp   i -1 , i   \rf\rp\rc.
		\end{align}
		By \cref{s:markovkoro} $M$ has a unique invariant distribution $\mu$ and since $N\stackrel{d}{=}\theta^{1}N$ it follows that $M_{0}\stackrel{d}{=}M_{1}\sim \mu$. The result follows from \cref{s:markovkoro} part 2.
		\\\\
\textbf{2.} Notice that $\alpha_{n}-\tau_{n}$ has $p$'th-moment /  exponential moment iff $\alpha_{n}-\lceil \tau_{n}\rceil$ has as well. For any realization such that $\tau_n< \infty$ we may write
		\begin{align}
		\alpha_{n} - \lceil\tau_{n}\rceil = \inf\lb i > 0 :\theta^{i-j+\lceil \tau_{n}-\alpha_{n-1}  \rceil}N_{\alpha_{n - 1} +}\lp  -1 , 0 \rf \leq \gamma\lp j\rp , j = 0,\dots,i +\lceil \tau_{n}-\alpha_{n-1}  \rceil- 1 \rb.
		\end{align}
		We now proceed as previously. Define  $M'_{i}$ for $i\geq 0 $ as the process
		\begin{align}
		M'_{i} = \inf\lb m \geq 0:\theta^{i-j+\lceil \tau_{n}-\alpha_{n-1}  \rceil}N_{\alpha_{n-1}+}\lp  -1, 0 \rf \leq \gamma\lp j+m\rp ,j = 0,\dots, i+\lceil \tau_{n}-\alpha_{n-1}  \rceil-1 \rb
		\end{align}
		 Notice that $\alpha_{n}-\lceil \tau_n\rceil$  is the first time $M'_{i}$ hits $0$. From \cref{rooflm} we have $\pi_{\alpha_{n-1}+}^{\uparrow \underline{\lambda},\overline{\lambda}} \indep \pi_{\alpha_{n-1}+}^{\downarrow \underline{\lambda},\overline{\lambda}}$ and thus $\tau_n - \alpha_{n-1}$ is independent of $N_{\alpha_{n-1}+}$. Observe also that $M'$ is an RE-process defined by
		\begin{align}
		M'_{i} = \lp M'_{i - 1}-1\rp \vee \lc\gamma^{ - 1}\lp N_{\tau_{n}+}\lp   i -1 , i  \rf\rp\rc .
		\end{align}
    To study the distribution of $M'_{0}$ note that it may be described as $M''_{\alpha_{n}-\lc \tau_{n}\rc}$ where $M''$ is another RE-process defined by $M''_{0}=0$ and
		\begin{align}
		M''_{i} = \lp M''_{i - 1}-1\rp \vee \lc \gamma^{ - 1}\lp N_{\alpha_{n-1}+}\lp   i -1 , i  \rf\rp\rc .
		\end{align}
		Let $P^*$ be the conditional distribution of $\lceil \tau_{n}-\alpha_{n-1}  \rceil$ given $\tau_n<\infty$. If $\phi$ is a positive increasing function and $\lp P^{k}_{x}\rp $ is the $k$-step Markov kernel for $M''$ then
		\begin{align}
		\qE\lp \phi\lp M'_{0}\rp \vert \tau_{n}<\infty  \rp   &= \int\int \phi\lp y\rp d P^{k}_{0}\lp y\rp dP^*\lp k\rp\\
		&\leq \int \phi\lp y\rp \mu\lp y\rp+\int\int \phi\lp y\rp d\lv P^{k}_{0}-\mu\rv\lp y\rp  dP^*\lp k\rp.
		\end{align}
		From theorem 14.1.4 \cite{MEYN} it follows that if $\mu\lp \phi\rp <\infty$, then the 2nd term above is finite. The desired result now follows from \cref{s:distributiontheorem} and \cref{s:markovkoro}.\\\\
\textbf{3.}		We now analyze $\alpha_{0}^I$ under setup (O). Here we utilize that the law of the stationary Linear Hawkes process $Z^0$ has a cluster process representation which we now describe; Let $N$ be a Poisson process on $\R$ with intensity $c_{\psi}$ and for $i\in\qZ$ let $(Z_{i})$ be independent Hawkes processes with weight/rate $h_{+},x\mapsto Lx$ and initialized with a single jump at $t=0$ (i.e. $Z_{i}\lf 0,0\rf=1$ and $R^{i}_{t}=h_{+}\lp t\rp $). Define now
		\begin{align}\label{eq:cluster}
		\overline{Z}(A) = \sum_{i\in \qZ\backslash \lb 0 \rb } C_{i}\lp A\rp
		\end{align}
where $C_{i}$ is a random measure given by $C_{i}\lp A\rp =Z_{i}\lp A-s_{i}\rp$ and $s_{i}$ is the i'th jump of $N$ for $i\in \qZ\backslash \lb 0 \rb $. Then $\overline{Z}$ is distributed as the stationary linear Hawkes process with weight/rate $h_{+}\slash \psi_{L}$. See \cite{DALEY},\cite{patroy} for more details on this construction. It follows that $\alpha^{I}_{0}\stackrel{D}{=}\tilde{\alpha}_{0}$ where
  \begin{align}
  \tilde{\alpha}_{0} &= \inf\{ i>0:\lv \int_{-\infty}^{i}h_{+}\lp t - s\rp  d\overline{Z}_{s} \rv  \leq f(t-i), \forall t>i \}.
  \end{align}
		Let $s_i$, $i\in \qZ\setminus\{0\}$ be the $i$th jump of $N$ before/after zero.  We shall use the following fact, coming from section 1.1 in \cite{patroy} and the proof of proposition 1.2 in same refrence: We may assume that the clusters $Z^i$ are constructed s.t. there is i.i.d $(W_{i},X_{i,j})$, also independent of $N$  such that $W_{i}\indep X_{i,j}$, $X_{i,j}\sim \lV h_{+} \rV^{-1}_{\cL^1} h_{+}(t)\ dt$ and
	\begin{align}\label{yw}
		Z_{i}(\R) = Z_{i}(-\infty,Y_{i}] = W_{i},
		\end{align}
		where $Y_{i} = \sum_{j=1}^{W_{i}}X_{i,j}$
and $W_{1}$ is distributed as the \textit{total progeny} of a Poisson branching process, with mean offspring $\lV h_{+} \rV^{-1}_{\cL^1}.$ \\\\ By the Otter-Dwass formula (see \cite{otterdwass}) The p.m.f. of $W_{1}$ is
\begin{align}
p_{W}\lp n\rp = \dfrac{n^{n-1}}{ \lV h_{+} \rV_{\cL^{1}} n!}e^{n\lp -\lV h_{+} \rV_{\cL^{1}}+\ln \lV h_{+} \rV_{\cL^{1}}\rp }.
\end{align}
 The stirling approximation for $n!$  gives that
\begin{align}\label{borelexp}
  \qE\exp\lp cW\rp<\infty \Longleftrightarrow c\leq c_{h}=\lV h_{+} \rV_{\cL^{1}}-\ln_{+} \lV h_{+} \rV_{\cL^{1}} -1
\end{align}
 If \cref{As:fmap} (B) holds, take any $c_{0}>1.$ Otherwise, let $c_{0}>1$ be a constant satisfying $\gamma\lp t\rp\geq c_{0}\lp p + 1\rp c_{h}^{ - 1}\ln_{ + }t$ for $t$ sufficiently large. Define $\gamma^{*}\lp t\rp=c_{0}^{-1}\gamma\lp c_{0}^{-1} t-1 \rp $ when $c_{0}^{-1}t \geq 1$ and $\gamma^{*}\lp t\rp = 0$ otherwise.\\\\
		With $N_{t}:=N[0,t]$ for $t>0$ and $-N[t,0]$ otherwise we define
		\begin{align}
		\overline{Y}_{i} &= \max_{l = N_{i-1}+1}^{N_i}\{Y_{l} \},\quad 		\overline{W}_{i} = \sum_{l = N_{i-1}+1}^{N_{i}}W_{l},\\
		\overline{\alpha}_{0} &= \inf\{i>0: \overline{Y}_{i-j} \leq \lp 1-c_{0}^{ - 1}\rp j,\ \overline{W}_{i - j} \leq \gamma^{*}\lp j\rp   ,\ j \geq 0\}.
		\end{align}
We claim that $\tilde{\alpha}_{0}\leq \overline{\alpha}_{0}$. This follows from the calculations for $t>\overline{\alpha}_0$
\begin{align}
	\int_{-\infty}^{\overline{\alpha}_0} h_+(t-s)d\overline{Z}_s
	&= \sum_{j=0}^{\infty} \int_{\overline{\alpha}_0-j-1}^{\overline{\alpha}_0-j} \int_{s}^t h_+(t-u) dC_{N_s}(u)dN_{s}\label{splitint}
\end{align}
 Note that \eqref{yw} implies that for all clusters $C_{N_{s}}$ with $s\in \lp \overline{\alpha_{0}}-j-1,\overline{\alpha_{0}}-j\rf$ we have  $Supp\lp C_{N_{s}}\rp \subset \lp s,s+\overline{Y}_{\overline{\alpha_{0}}-j}\rf$. By the inequalities obtained from the definition of $\overline{\alpha}_0$ we have $\lp s,s+\overline{Y}_{\overline{\alpha_{0}}-j}\rf\subset \lp s,\overline{\alpha_{0}}-c_{0}^{-1}j\rf.$  Since $h_{+}\leq \overline{h}$ and $\overline{h}$ is decreasing we conclude that $\overline{h}\lp t-u\rp \leq \overline{h}\lp t-\overline{\alpha}_0 + c_0^{-1}j\rp $ for all $u\in Supp\lp C_{N_{s}}\rp,s\in \lp \overline{\alpha_{0}}-j-1,\overline{\alpha_{0}}-j\rf,t>\alpha_{0}$. Inserting this into  \eqref{splitint} gives
\begin{align}
\int_{-\infty}^{\overline{\alpha}_0} h_+(t-s)d\overline{Z}_s
&\leq \sum_{j=0}^{\infty} \overline{h}(t-\overline{\alpha}_0 + c_0^{-1}j)  \int_{\overline{\alpha}_0-j-1}^{\overline{\alpha}_0-j} C_{N_s}[s,\overline{\alpha}_0 - c_0^{-1}j] dN_s.\\
&\leq \sum_{j=0}^{\infty}\overline{h}(t-\overline{\alpha}_0 + c_0^{-1}j) \overline{W}_{\overline{\alpha}_0-j}
\end{align}
Again, by the inequalities defining $\overline{\alpha}_{0}$ we get
\begin{align}
\int_{-\infty}^{\overline{\alpha}_0} h_+(t-s)d\overline{Z}_s &\leq \int_{0}^{\infty} \overline{h}(t-\overline{\alpha}_0 + c_0^{-1}s) \gamma^*(s+1) ds + \overline{h}(t-\overline{\alpha}_0)\gamma^*(0)\\
&\leq \int_{0}^{\infty} \overline{h}(t-\overline{\alpha}_0 + s) \gamma(s+1) ds + \overline{h}(t-\overline{\alpha}_0)\gamma(0)\\
&\leq f(t-\overline{\alpha}_0).
\end{align}
which proves the claim.

To describe the moment of $ \overline{\alpha}_{0}$ define
	\begin{align}
	M_{i} = \inf\lb m \geq 0:\overline{Y}_{i-j} \leq \lp 1-c_{0}^{ - 1}\rp \lp j+m\rp ,\ \overline{W}_{i - j} \leq \gamma^{*}\lp j+m\rp  ,j \geq 0 \rb \label{MinvZ}.
	\end{align}
		Indeed, $M_{0}<\infty$ by a Borel-Cantelli argument. As before $\overline{\alpha}_{0}$ is the return time to $0$ for the RE-process  $M_{i}$ with update-variables
\begin{align*}
  \lc\lp \lp 1-c_{0}^{-1}\rp^{-1}  \overline{Y}_i\rp \vee \lp \gamma^{*}\rp^{-1} \lp \overline{W}_{i}\rp\rc
\end{align*}
 and started at $M_{0}$, carrying the invariant distribution of $M$. Under \cref{As:fmap} $(A)$ we have $\lp \gamma^*\rp^{-1}\lp t\rp\leq \exp \lp c_{h}\lp p+1\rp^{-1}t\rp$ for $t$ large. It follows that the update variables have $(p+1)$th moment and the starting distribution has $p$'th moment, so it follows from \cref{s:markovkoro} that $\overline{\alpha}_{0}$ has $p$'th moment. Likewise, under \cref{As:fmap} $(B)$ the update variables have exponential moment and the starting distribution has exponential moment as well. It follows from \cref{s:markovkoro} that $\overline{\alpha_{0}}$ has exponential moment.\\\\
\\\textbf{4.}		We now analyze $\alpha_{n}$ under setup (O). To outline the similarity we shall re-use some of the notation from the previous proof, but for slightly modified random variables.

		Notice that $Z^n[t,t]$ is dominated by the Hawkes process  $Z'$ with the same weight function $h_+$, rate function $\psi_{L}$  and initial signal $R'_{t}=f(t-\lff t\rff) + h_{+}(t-\lp \tau_{n}-\alpha_{n-1}\rp )$.  The law of $Z'$ has a cluster process representation given as follows: Let $N$ be an inhomogeneous Poisson Process with intensity $c_{\psi}+Lf(t-\lff t\rff)$, and let $\xi\sim (\tau_n-\alpha_{n-1})\mid (\tau_n<\infty)$. Define $Z_{i}$ as before, and let $C_i$, $C^{\xi}$ be the mutually independent clusters given by $C_i(A)=Z_{i}(A-s_i)$ where $s_i$ is the $i$th jump of $N$ for $i\in \N$, and $C^{\xi}(A)=Z_{-1}(A-\xi)$. Define
		\begin{align}
			\overline{Z}(A) = \sum_{i=0}^{\infty} C_i (A) + C^{\xi}(A), \quad A\in \cB_{\R_+}.
		\end{align}
		Then $\overline{Z}\stackrel{D}{=}Z'$ and hence also
{\small
		\begin{align}
			&\lp \alpha_n-\alpha_{n-1}\mid \tau_n<\infty \rp\\ \stackrel{D}{=} &\inf \left\{i> \lc \xi \rc  :\lv \int_{0}^{i}h_{+}\lp t - s\rp  d\overline{Z}_{s}  \rv  \leq f\lp t - i , i-1 \rp, \forall t>i,\overline{Z}\lf 0,i\rf  \leq \int_{0}^{i} \gamma\lp s+1\rp ds  \right\}. \label{zeta}
		\end{align}
}
  as before, we have i.i.d. variables
		$(W_{i},X_{i,j}),(W^{\xi}, X^{\xi}_j)$ and also independent of $N, \xi$  such that $X_{i,j}\sim \lV h_{+} \rV^{-1}_{\cL^1} h(t)\ dt$, $W_{i}\sim Z_{0}(\R)$ and
		\begin{align}
		Z_{i}(\R) &= Z_{i}(-\infty,Y_i] = W_{i}, \quad \forall i \in \N\\
		Z_{-1}(\R) &= Z_{-1}(-\infty, Y^{\xi}] = W^{\xi}
		\end{align}
		where $Y_i= \sum_{j=1}^{W_{i}}X_{i,j}$ and $Y^{\xi}=\sum_{j=1}^{W^{\xi}}X^{\xi}_j$. Define now
		\begin{align}
			\overline{Y}_i = \max_{l=N_{i-1}+1}^{N_i} \{Y_l\}, \quad \overline{W}_i = \sum_{l=N_{i-1}+1}^{N_i} W_{l}
		\end{align}
		and notice that $(\overline{Y}_i, \overline{W}_i)$ for $i\in\N$ is an i.i.d. sequence. Define now $c_{0},\gamma^*$ as previously and set
		\begin{align}
			\overline{\zeta}
			= \inf\{i>0 :\  &\overline{Y}_{i+\lceil \xi\rceil-j} \leq \lp 1-c_{0}^{ - 1}\rp j,\ \overline{W}_{i+\lceil \xi\rceil - j} \leq \gamma^{*}\lp j\rp \ \forall j \in \{0, ...,i+\lceil\xi\rceil-1\}\setminus\{i\}, \\
			&\overline{Y}_{\lceil \xi\rceil}\vee Y^{\xi} \leq \lp 1-c_{0}^{ - 1}\rp i,\ \ \overline{W}_{\lc \xi\rc} + W^{\xi}\leq \gamma^{*}\lp i\rp\}.
		\end{align}
		Define now the random time $\zeta$ s.t. $\zeta+\lc \xi\rc $ is equal to \eqref{zeta}. We claim that $\zeta  \leq  \overline{\zeta}$. To prove this, we apply similar arguments as in the part of \textbf{3}) where we showed that \\ $\tilde{\alpha}_{0}\leq \overline{\alpha}_{0}.$ By construction of $\overline{Y}_{i}$'s and $\overline{\zeta}$, it holds for all $s\in \lp \overline{\zeta}+\lf \xi \rf-j-1,\overline{\zeta}+\lc \xi \rc-j\rp $ that  $Supp \lp C_{N_{s}}\rp\subset \lp -\infty, \overline{\zeta}+\lc \xi \rc-c_{0}^{-1}j \rf$ and this implies
		\begin{align}
		&\int_{0}^{\lc\xi\rc + \overline{\zeta}} h_+(t-s)d\overline{Z}_s\\
		\leq &\sum_{j=0,j\neq \overline{\zeta}}^{\lc\xi\rc + \overline{\zeta} - 1}\overline{h}(t-\lp \lc\xi\rc + \overline{\zeta}\rp + c_0^{-1}j) \overline{W}_{\lc\xi\rc + \overline{\zeta}-j} + \overline{h}(t-\lp \lc\xi\rc + \overline{\zeta}\rp  + c_0^{-1}\overline{\zeta})\lp \overline{W}_{\lc\xi\rc} + W^{\xi}\rp \\
		\leq &\int_{0}^{\lc\xi\rc+\overline{\zeta}-1} \overline{h}(t-\lp \lc\xi\rc + \overline{\zeta}\rp  + c_0^{-1}s) \gamma^*(s+1) ds + \overline{h}(t-\lp \lc\xi\rc + \overline{\zeta}\rp )\gamma^*(0)\\
		\leq &f(t-\lp \lc\xi\rc + \overline{\zeta}\rp ,\lc\xi\rc + \overline{\zeta}-1).
		\end{align}
Likewise, the support constraint on $C_{N_{s}}$ mentioned above must imply that
\begin{align}
  \overline{Z}\lp 0, \overline{\zeta}+\lc \xi\rc \rf = \sum_{j=0}^{\overline{\zeta}+\lf \xi \rf-1} W_{j} \leq \sum_{j=0}^{\overline{\zeta}+\lc \xi\rc-1}\gamma^*\lp j\rp\leq \int_{0}^{\overline{\zeta}+\lc \xi\rc-1} \gamma\lp s+1\rp ds.
\end{align}
This proves our claim.
		To analyze $\overline{\zeta}$ define
		\begin{align}
			M_i = \inf\lb m \geq 0:\overline{Y}_{i-j} \leq \lp 1-c_{0}^{ - 1}\rp \lp j+m\rp ,\ \overline{W}_{i - j} \leq \gamma^{*}\lp j+m\rp ,\ j = 0, ..., i-1 \rb.
		\end{align}
		This is an RE-process started at $M_0=0$ and with update variables \\
  $\lc \max\{(\gamma^*)^{-1}(\overline{W}_{i}), (1-c_0^{-1})^{-1}\overline{Y}_{i}\}\rc $. As in the proof of \textbf{2)} it follows that $M_{\lf \xi\rf -1}$ has $p$'th moment, and under assumption 2 B) it has exponential moment. The same therefore holds for the variable
		\begin{align}
			M'_0 :=  (M_{\lf \xi\rf -1} - 1) \vee \lc \max\{(\gamma^*)^{-1}(\overline{W}_{\lf \xi\rf }+W^{\xi}), (1-c_0^{-1})^{-1}\overline{Y}_{\lf \xi\rf } \vee Y^{\xi} \}\rc .
		\end{align}
		Now realize that $\overline{\zeta}$ is the return time to 0 for the RE-process with update variables
\begin{align}
  \lc \left(\left(1-c_0^{-1}\right)^{-1}\overline{Y}_i\right) \vee (\gamma^*)^{-1}(\overline{W}_i)\rc
\end{align}
 and started at $M'_0$. The desired result now follows from \cref{s:markovkoro}
	\qed

\vspace{1cm}
\noindent\textbf{Proof of \cref{regenerationTheorem}}\\

\textbf{1.} To prove that  $\rho$ is a $\cF^*$ stopping time, we notice that
	\begin{align}
	(\rho\leq t)=(\alpha_{\eta}\leq t-D)= \bigcup_{i=0}^{\lfloor t-D \rfloor} (\alpha_{\eta} = i).
	\end{align}

	So it suffices to show $(\alpha_{\eta}=i) \in \cF^*_i$. Indeed, this is true since
	\begin{align}
	(\alpha_{\eta}=i) = \bigcup_{k=0}^{\infty} (\alpha_k = i) \cap \lp \int_{i}^{\infty} \ql \{z\leq F(s)\}d\pi^{\uparrow \underline{\lambda},\overline{\lambda}} \rp.
	\end{align}

	\textbf{2.} By applying \cref{rooflm} for each fixed $t$ we see that  $\pi^{\downarrow \underline{\lambda},\overline{\lambda}}$ is an $\cF^*_{t}$-PRM. The Strong Markov Property gives that $\pi_{\rho+}^{\downarrow \underline{\lambda},\overline{\lambda}}$ is a PRM independent of $\cF^*_{\rho}$. Notice now that $Z_+,\lambda_{+t}:=\lambda_{\rho+t}$ is exactly the point process and intensity of the ADHP driven by $\pi_{\rho+}^{\downarrow \underline{\lambda},\overline{\lambda}}$ with initial age $D$ and signal $R:t\mapsto -f\lp t+D\rp $. That is, $Z_{\rho+}$ is entirely generated by $\pi_{\rho+}^{\downarrow \underline{\lambda},\overline{\lambda}}$ and hence $Z_{\rho+}\indep \cF^*_{\rho}$. Notice that $Z_{\rho+}=Z_{\rho+}^*$ by \cref{s:bandlemma} and that $Z^*_{\vert \lp 0,\rho \rf}\subset \cF^*_{\rho}$. It now follows that $Z^*_{\rho+}\indep Z^*_{\vert \lp 0,\rho \rf}.$
 \\

\textbf{3.} We introduce an i.i.d. sequence of random variables  $\left( \beta_{i}\right)$ with distribution \\ $\beta_{1}\sim \alpha_{1} - \alpha_{0}\vert \tau_{1}-\alpha_{0}<\infty$. We then introduce another sequence of random variables given by
	\begin{equation}{
		\tilde{\beta}_{i} =
		\begin{cases}
		\alpha_{i} - \alpha_{i - 1}& i\leq \eta  \\
		\beta_{i}& i>\eta.
		\end{cases}
	}\end{equation}
	Recall that $\lp \alpha_{j} - \alpha_{j - 1}\rp_{1\leq j\leq i}$ are conditionally i.i.d. given $\tau_{i}<\infty.$ From this we see that $\tilde{\beta}$ is an i.i.d. sequence of variables distributed as $\beta_{1}$, and which is also independent of $\eta $. We may now write
	\begin{align}{
		\alpha_{\eta}-\alpha_{0}&= \sum_{i = 1}^{\eta} \alpha_{i} - \alpha_{i  - 1}\\
		&=\sum_{i = 1}^{\eta} \tilde{\beta}_{i}.
	}\end{align}
	From  \cref{s:distributiontheorem} 1) it follows that $\eta$ is distributed as a negative binomial, and in particular it has exponential moment. To study the distribution of $\tilde{\beta}_{1} $, we use that
	\begin{align}
	\alpha_{1} - \alpha_{0} = \lp \alpha_{1} - \tau_{1}\rp  + \lp \tau_{1} - \alpha_{0}\rp.
	\end{align}
 	 \Cref{s:distributiontheorem} and \cref{s:AlphaDist} gives that $\beta$ has p'th moment so the desired result follows from theorem 5.2, chapter 1 \cite{Gut}. Under Assumption 2B, one may conclude the proof by writing  for small $c>0$
 	\begin{align}
 		\qE \exp(c \alpha_{\eta})
 		&= \qE \exp(c\alpha_{0}) \ \qE \prod_{i=0}^{\eta} \exp(c\tilde{\beta}_i) \\
 		&= \qE \exp(c\alpha_{0}) \ \qE \left[\qE\left( \prod_{i=0}^{\eta} \exp(c\tilde{\beta}_i) \ \big\vert\ \eta \right)\right] \\
 		&=\qE \exp(c\alpha_{0}) \ \qE \left(\qE\exp(c\tilde{\beta}_1)\right)^{\eta}.
 	\end{align}
 	Since $\eta$ has exponential moment, the above expression is finite for small $c$. \qed

\subsection{Proofs of Section 3.2 Results}
\noindent\textbf{Proof of \cref{prop:inv}}\\
The proof is a coupling argument. Consider the Hawkes process $Z^*$ driven by a fixed PRM $\pi$ and started with $R:t\mapsto f\lp t+D\rp ,A_{0}=D$.  Then \cref{As:alpha0} is satisfied with $\alpha_{0}=0$, $ r=-f.$ Thus, defining $\Phi$ as in \eqref{eq:markovlol} we have a coupling in the sense that $\lp \theta^{n} Z^*\rp_{\vert \lp -D,0\rf}=H^*\lp \Phi_{n}\rp$ for all $n>\rho_{0}$. Since $\tilde{P}$ is the invariant distribution of $\Phi$ the markov chain converges in total variation and hence also $H^*\lp \Phi_{n} \rp\Rightarrow H^*\lp \tilde{P}\rp.$ On the other hand, we have by the coupling property of $Z^I$ that almost surely, there is a random integer $n_{0}$ s.t. \\ $\lp \theta^n Z^*\rp_{\vert \lp -D,0\rf}=\lp \theta^n Z^I\rp_{\vert \lp -D,0\rf}$ for all $n\geq n_{0}$. It follows that $H^*\lp \Phi^n \rp=\lp \theta^n Z^*\rp_{\lp -D,0\rf}\Rightarrow Z^I_{\lp -D,0\rf}$ and the desired result follows from uniqueness of limits.
\qed


\section{Appendix}
\subsection{Splitting Two PRM's}
Let  $\pi,\overline{\pi}$ be two independent PRM's on $\R^2_{+}$. For any two functions \\$ f_1: \R \con \lf 0,\infty \rp , f_{2}:\R\con \lf 0,\infty \rf$, such that  $f_{1}\leq f_{2}$ we define for $B\in \cB_{\R^2_{+}}$
\begin{align}
\pi^{\downarrow f_{1},f_{2}}\left( B\right)  &=\int_{B}\ql\left\{ z \not\in \lf f_{1},f_{2} \rp   \right\} \;d\pi\left( s,z\right) + \int_{B}\ql\left\{ z \in  \lf f_{1},f_{2}\rp   \right\} \;d\overline{\pi}\left( s,z\right).\\
\pi^{\uparrow f_{1},f_{2}}\left( B\right)
&=\int\ql\left\{ \lp s,z\rp : \lp s,z-f_{1}\lp s\rp\rp  \in B,\ z\leq f_2(s) \right\} \;d\pi\left( s,z\right)\\
&\hspace{4mm}+ \int\ql\left\{ \lp s,z\rp : \lp s,z-f_{1}\lp s\rp\rp  \in B,\ z> f_2(s) \right\} \;d\overline{\pi}\left( s,z\right).
\end{align}
It may be shown directly that both of the above set functions are PRM's with $\pi^{\downarrow f_{1},f_{2}}\indep \pi^{\uparrow f_{1},f_{2}}$.

\Cref{rooflm} shows that these independence properties generalizes to when $f_{1}, f_2 $ are predictable, but random, intensities.

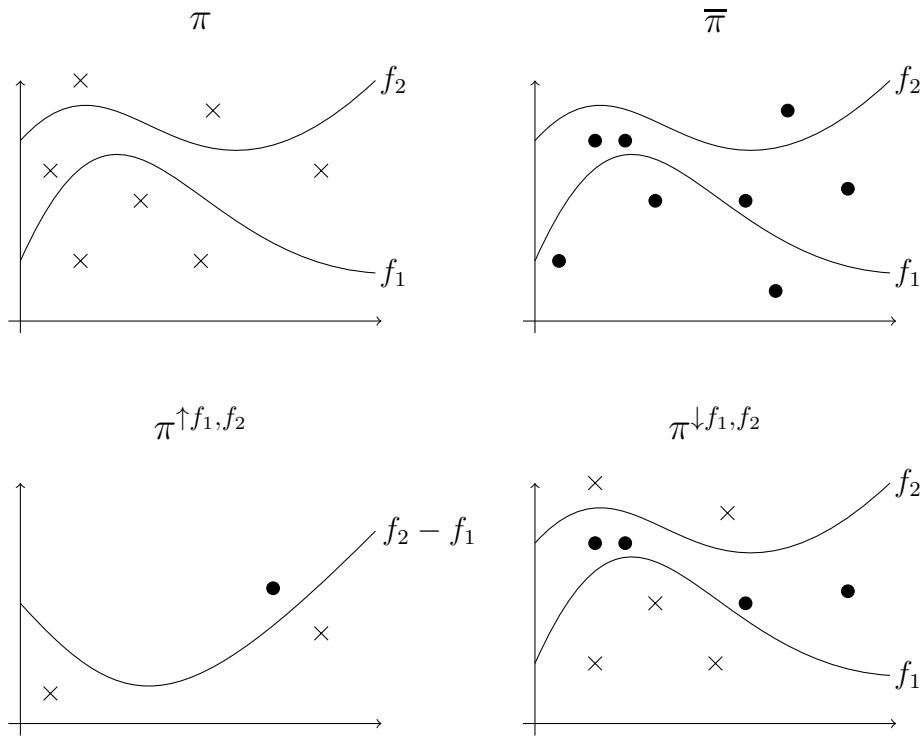
\begin{figure}[H]
	\begin{center}
		\begin{tikzpicture}[scale=0.8]
		\draw (3,5) node {\large$\pi$};
		\draw[->] (-0.2,0) -- (6,0);
		\draw[->] (0,-0.2) -- (0,4);
		\draw (8, 0) node {};

		\draw (0,3) ..controls (1.8,5) and (2.7, 1) .. (5.9,4);
		\draw (6.2,4) node {$f_2$};

		\draw (0,1) ..controls (1.8,5) and (2.7, 1) .. (5.9,0.8);
		\draw (6.2,0.8) node {$f_1$};

		\draw (1,1) node[cross] {};
		\draw (0.5,2.5) node [cross] {};
		\draw (1,4) node[cross] {};
		\draw (3.2,3.5) node[cross] {};
		\draw (2,2) node[cross] {};
		\draw (3,1) node[cross] {};
		\draw (5,2.5) node[cross] {};
		\end{tikzpicture}
		\begin{tikzpicture}[scale=0.8]
		\draw (3,5) node {\large$\overline{\pi}$};
		\draw[->] (-0.2,0) -- (6,0);
		\draw[->] (0,-0.2) -- (0,4);

		\draw (0,3) ..controls (1.8,5) and (2.7, 1) .. (5.9,4);
		\draw (6.2,4) node {$f_2$};

		\draw (0,1) ..controls (1.8,5) and (2.7, 1) .. (5.9,0.8);
		\draw (6.2,0.8) node {$f_1$};

		\filldraw (0.4,1) circle (3pt);
		\filldraw (1,3) circle (3pt);
		\filldraw (2,2) circle (3pt);
		\filldraw (1.5,3) circle (3pt);
		\filldraw (3.5,2) circle (3pt);
		\filldraw (4.2,3.5) circle (3pt);
		\filldraw (5.2,2.2) circle (3pt);
		\filldraw (4,0.5) circle (3pt);
		\end{tikzpicture}

		\vspace{0.8cm}

		\begin{tikzpicture}[scale=0.8]
		\draw (3,5) node {\large$\pi^{\uparrow f_1, f_2}$};
		\draw[->] (-0.2,0) -- (6,0);
		\draw[->] (0,-0.2) -- (0,4);
		\draw (8, 0) node {};

		\draw (0,2) ..controls (1.8,0) and (2.7, 0) .. (5.9,3.2);
		\draw (6.8, 3.2) node {$f_2-f_1$};
		\draw (7, 0) node {};

		\draw (0.5,0.5) node [cross] {};
		\draw (5,1.5) node[cross] {};
		\filldraw (4.2,2.25) circle (3pt);

		\end{tikzpicture}
		\begin{tikzpicture}[scale=0.8]
		\draw (3,5) node {\large$\pi^{\downarrow f_1, f_2}$};
		\draw[->] (-0.2,0) -- (6,0);
		\draw[->] (0,-0.2) -- (0,4);

		\draw (0,3) ..controls (1.8,5) and (2.7, 1) .. (5.9,4);
		\draw (6.2,4) node {$f_2$};

		\draw (0,1) ..controls (1.8,5) and (2.7, 1) .. (5.9,0.8);
		\draw (6.2,0.8) node {$f_1$};

		\draw (1,1) node[cross] {};
		\draw (2,2) node[cross] {};
		\draw (3,1) node[cross] {};
		\draw (1,4) node[cross] {};
		\draw (3.2,3.5) node[cross] {};
		\filldraw (1,3) circle (3pt);
		\filldraw (1.5,3) circle (3pt);
		\filldraw (3.5,2) circle (3pt);
		\filldraw (5.2,2.2) circle (3pt);
		\end{tikzpicture}

	\end{center}
\caption{ \label{fig:pi} A figure illustrating how $\pi^{\uparrow f_1, f_2}$ and $\pi^{\downarrow f_1, f_2}$ are created from $\pi, \overline{\pi}$. Notice that $\pi^{\downarrow f_1, f_2}$ contains the part of $\pi$ \textit{below} $f_{1}$ while $\pi^{\uparrow f_1, f_2}$ contains an area part immediately \textit{above} it.}
\end{figure}

\begin{s}\label{rooflm}
	Let $(\cF_t)_{t\in\R_+}$ be a filtration, and $\pi, \overline{\pi}$ be two independent $\cF_t$-PRMs on $\R_{+}$. For $t\in \lf 0,\infty \rp$, let $\lambda_{t}\leq \lambda_{t}'$ be $\cF_{t}$-predictable processes taking values in $\lf 0,\infty\rp,\lf 0,\infty\rf$, respectively. Define $\cF^*_{t}=\sigma\lp \cF_{t},\pi^{\uparrow \lambda, \lambda'}\rp $ It holds that $\pi^{\uparrow \lambda, \lambda'}, \pi^{\downarrow \lambda, \lambda'}$ are PRM's such that  $\pi^{\uparrow \lambda, \lambda'}\indep \pi^{\downarrow \lambda, \lambda'}\indep \cF_{0}.$
\end{s}

\begin{p}
	The proof will be done in several steps.\\

	\noindent \textbf{Step 1.}\\
	Let $\lp  t_{i}\rp_{i \in \N_{0}} $ be a fixed partition  where $0<t_{i-1} < t_i $. $t_i\conx{i\con\infty}\infty$. Assume $ Y^1_i,Y^2_i $ is ${\mathcal F}_{t_{i-1}}-$measurable, taking values in a finite state space $\mathcal{Y}\subset \R$, and assume that

	\begin{align}
	\lambda_t = \sum_{i=1}^\infty Y^1_i 1_{]t_{i-1}, t_i ] }(t).\quad 	\lambda'_t = \sum_{i=1}^\infty Y^2_i 1_{]t_{i-1}, t_i ] }(t). 
	\end{align}

	Fix $m\in\N$ and take $k^{i}_{j},l^i_j \in \N$ and mutually disjoint, bounded $B^{i}_{j}, C^i_j \in \mathcal{B}_{\left( t_{i - 1},t_{i} \times\R_{+} \right]}$ for $i = 1,\dots,n\;\;j = 1,\dots,m$.
	Define
	\begin{align}
	F_{i} = \bigcap_{j = 1}^{m}\left( \pi^{\downarrow\lambda,\lambda'}\left( B^{ i}_{j}\right) = k^{i}_{j}\right), \\
	G_{i} = \bigcap_{j = 1}^{m}\left( \pi^{\uparrow\lambda,\lambda'}\left( C^{ i}_{j}\right) = l^{i}_{j}\right).
	\end{align}
Take $E\in \cF_{0}.$ and set $E_{i} = F_{i}\cap G_{i}, E=\bigcap_{i=0}^n E_{i}.$ It is sufficient to show that the projection $E$ has the correct distribution, i.e.:
	\begin{align}
		P\lp E\rp = P\lp E_{0}\rp \prod_{i = 1}^{n}\prod_{j = 1}^{m}\mathcal{P}\left( \int_{B^{i}_{j}}ds,k^{i}_{j}\right)\mathcal{P}\left( \int_{C^{i}_{j}}ds,l^{i}_{j}\right)
	\end{align}
	where $\mathcal{P}\left( c,\cdot\right)  $ is the Poisson density with mean $c $. This will be proved using induction. The induction claim over $N=0,...,n$ is that
	\begin{align}{
		P\lp E_{0}\rp   = \mathbb{E}\ql\left\{ E_{0}\cap  E_{1}\cap  \dots\cap E_{n-N}\right\}\prod_{i=n-N+1}^n \prod_{j = 1}^{m}\mathcal{P}\left( \int_{B^{i}_{j}}ds,k^{i}_{j}\right)\mathcal{P}\left( \int_{C^{i}_{j}}ds,l^{i}_{j}\right).\label{condprop}
	}\end{align}
The induction start $N=0$ is clear (where the empty product is 1 per convention). Assume that the claim holds for some $N-1$.
	Since  $E_{N}\indep \cF_{t_{N - 1}}\vert \lambda_{t_{N}},\lambda'_{t_{N}}$
	we may write \\ $P\left( E_{N}\vert \cF_{t_{N - 1}}\right) = P\left( E_{N}\vert \lambda_{t_{N}},\lambda'_{t_{N}}\right)$ and for $c\leq d \in \mathcal{Y} $
	\begin{align}{
		&P\left( E_{N}\vert \lp \lambda_{t_{N}},\lambda'_{t_{N}}\rp =\lp c,d\rp \right)\\
		&= P\bigcap_{j = 1}^{m}\left( \int_{B_{j}^{N}}\ql\left\{ z \not\in \lf c,d\rp   \right\} \;d\pi\left( s,z\right) + \int_{B_{j}^{N}}\ql\left\{ z \in  \lf c,d \rp   \right\} \;d\overline{\pi}\left( s,z\right)=k^{i}_{j} \right)\\
		&\cdot P\bigcap_{j = 1}^{m}\left( \int_{c+C_{j}^{N}}\ql\left\{ z \in \lf c,d\rp   \right\} \;d\pi\left( s,z\right) + \int_{c+C_{j}^{N}}\ql\left\{ z \not\in  \lf c,d \rp   \right\} \;d\overline{\pi}\left( s,z\right)=l^{i}_{j} \right).\label{prodpois}
	}\end{align}
Straightforward calculations show that indeed

	\begin{align}{
		P\left( E_{N}\vert  \lp \lambda_{t_{N}},\lambda'_{t_{n}}\rp =\lp c,d\rp\right) = \prod_{j=1}^{m}\mathcal{P}\left( \int_{B^{N}_{j}}ds,k^{N}_{j}\right)\mathcal{P}\left( \int_{C^{N}_{j}}ds,l^{N}_{j}\right).
	}\end{align}
	Notice that the right side does not depend on $c,d$ implying that $E_{N}\indep \cF_{t_{N-1}}$.Thus, by conditioning \eqref{condprop} w.r.t. $\cF_{t_{N-1}}$ and inserting this result, the induction step follows and the proof is completed.\\

	\noindent\textbf{Step 2.}\\
	Assume now that $\lambda,\lambda' $ is bounded and continuous in $t$ for all discrete measures. We use some dyadic approximation by putting
	$$  \lfloor x\rfloor_n := \sup_{k} \{ k 2^{-n} : k 2^{-n } < x \}  $$
	$$ \lambda^{n}_{t} = \lfloor \lambda_{\lfloor t\rfloor_n} \rfloor_n,\quad \lambda^{'n}_{t} = \lfloor \lambda'_{\lfloor t\rfloor_n} \rfloor_n.$$

	Then   as $ n \to \infty ,$

	\begin{align}{
		\lambda^{n}_{t}\con \lambda_{t}, \quad \lambda^{'n}_{t}\con \lambda'_{t}\label{Hconverge}
	}\end{align}
	for all $t \in \left[ 0,T\right].$
	Almost surely, the graphs of $\lambda,\lambda'$ are $\pi,\overline{\pi}$ null-sets. It follows that
	almost surely

	\begin{align}
	\forall A\in \cB^{2},Leb\lp A\rp<\infty :\pi^{\uparrow\lambda^{n},\lambda^{'n}}\lp A\rp  \con \pi^{\uparrow\lambda,\lambda^{'}}\lp A\rp ,\quad \pi^{\downarrow\lambda^{n},\lambda^{'n}}\lp A\rp  \con \pi^{\downarrow\lambda,\lambda^{'}}\lp A\rp.
	\end{align}
	It is now straightforward to prove the claim by applying step 1 for each $n$.\\

	\noindent\textbf{Step 3.}\\
	Assume the same set up as before, except for continuity in $t $. Define for $n \in \N $

	\begin{align}
	\lambda^{n}_t = n\int_{t-\frac{1}{n}}^{t} \lambda_s ds,
	\quad \lambda'^n_t =n\int_{t-\frac{1}{n}}^{t} \lambda'_s ds.
	\end{align}

	Note that the processes

	\begin{align}
	t \mapsto \int_{0}^{t}\lambda_s ds, \quad t \mapsto \int_{0}^{t} \lambda'_s ds
	\end{align}
	is Lipschitz continuous since $\lambda, \lambda'$ is bounded. By Rademacher's theorem there is a Lebesque-full set on which the above map is differentiable. It follows that almost surely, $\lambda^n, \lambda'^n$ converges almost everywhere in $t$ to $\lambda, \lambda'$. The remaining part of step $2 $ is similar to step $3 $.\\

	\noindent\textbf{Step 4.}\\
	Assume now that $\lambda, \lambda'$ are given as in the assumptions. One may define $\lambda_n=\lambda \wedge n,\lambda'_n=\lambda'\wedge n$, and repeat the procedure from the previous steps to complete the proof, which we leave to the reader.
\end{p}

\subsection{The Random Exchange Process}
The purpose of this section is to study the Markov Chain given by
\begin{align}
M_{i} = \lp M_{i-1} - 1\rp \vee X_{i}
\end{align}
where $M_{0},X_{i}$ are non-negative and mutually independent variables such that $\lp X_{i}\rp $ are i.i.d. This process is going under the name Random Exchange Process with constant decrements. RE-processes have been treated in \cite{sill} where it was shown that $M$ is positive recurrent when $X$ has finite expectation. See also \cite{null} for a null-recurrence characterization. We are interested in moments of the return time $\sigma=\inf \lb n>0:M_{n}\leq 0 \rb $, and moments of the invariant distribution $\mu$. To the best knowledge of this author, there has been no published result about such.\\\\
Let $F,S$ be the distribution function and survival function of $X_{1}$. Let $q\geq 0$ be a real number. Clearly the transition kernel of $M_{i}$ is given by
\begin{align*}
P_{x}\lp \lp a,b\rf \rp &= P\lp X\in \lp a,b\rf \rp\quad  \forall  b>a>x-1  \\
P_{x}\lp \lb x-1\rb  \rp &= F\lp x-1\rp.
\end{align*}
 Let $\phi: \lf 0,\infty\rp \con \lf 0,\infty\rp$ denote an increasing and differentiable function. Valid choices of $\phi$ include $\phi\lp x\rp=x^{q+1}$ and $\phi\lp x\rp=\exp\lp cx\rp.$
\begin{s}\label{s:lyapounov}
Assume that $\phi$ is convex and $\qE \phi \lp X\rp<\infty.$
 Then $M_{i}$ is positive recurrent with stationary distribution $\mu $ and $\int \phi'\lp y-1\rp d\mu\lp y\rp<\infty$.
\end{s}
\begin{p}
We use a Lyapounov argument.
\begin{align}
\qE_{x} \phi\lp M_{1}\rp - \phi\lp x\rp  &= \int_{\lp x-1,\infty\rp} \phi \lp y\rp dF\lp y\rp   + \phi\lp x - 1\rp F\lp x-1\rp  - \phi\lp x\rp\\
& = \int_{\lp x-1,\infty\rp} \phi \lp y\rp dF\lp y\rp + \lp \phi\lp x - 1\rp - \phi\lp x\rp \rp  F\lp x-1\rp - S\lp x-1\rp \phi\lp x\rp.
\end{align}
Notice that the first term above converges to $0$ for $x\con \infty $. The second term can be controlled using the mean value theorem. Indeed, for $x $ so large that $F\lp x-1 \rp \geq 2^{ - 1}  $ we have

\begin{align}
\lp \phi\lp x - 1\rp - \phi\lp x\rp \rp  F\lp x-1\rp  \leq  -  \dfrac{1}{2}\phi'\lp x - 1\rp
\end{align}
We may apply  Proposition 14.1.1 and theroem 14.2.3 from \cite{MEYN} with $f\lp m\rp = \dfrac{1}{2}\phi'\lp x - 1\rp $ to conclude the desired result. \end{p}

\noindent Note that $\mu$ must satisfy $\mu\lf 0,x\rf =F\lp x\rp \mu\lf 0,x+1\rf$, and in turn it satisfies
\begin{align}
\mu\lp \lf 0,x \rf \rp  = \prod_{k=0}^{\infty}F\lp x+k\rp
\end{align}
whenever it exists. In particular when $X$ has support on $\N_{0}$ we have
\begin{align}
\mu\lp \lf 0,n \rf \rp  = \prod_{k=n}^{\infty}F\lp k\rp.
\end{align}
We now discuss the hitting time $\sigma$. We need an intermediate result that gives a peculiar relation between the return time to $0$, and the hitting time given general distributions for a RE-process,  when the update variables have support on $\N_{0}$.
\begin{s}\label{s:hittingtimesum}
Assume that $X$ has support on $\N_{0}.$ Define for $i,j\in\N_{0}$ $e_{i,j}=\qE_{i}\phi\lp \sigma+j\rp$. Assume that $e_{i,j}<\infty$  for all $i,j$.
 For any probability measure $\nu$ on $\N_{0}$ it holds that
\begin{align}
\qE_{\nu}\phi\lp \sigma\rp=
\phi\lp 0\rp + \nu\lp 0\rp \lp \qE_{0}\phi\lp \sigma\rp-\phi\lp 0\rp \rp + \sum_{i=0}^{\infty} \left(\qE_{1}\phi(\sigma + i) - \phi(i)\right)S_{\nu}(i)\frac{\mu\lp 0\rp }{\mu\lf 0, i\rf }
\end{align}
where $S_{\nu}$ is the survival function of $\nu$.
\end{s}
\begin{p}
note that
	\begin{align}
	\qE_{\nu}\phi\lp \sigma\rp=\sum_{i=1}^{\infty}\nu\lp i\rp e_{i0}. \label{eq:expectationformula}
	\end{align}
	We claim that
	\begin{align}
	e_{i,j}-e_{i-1,j}=F\lp i-2\rp \lp  e_{i-1,j+1}-e_{i-2,j+1}\rp , \quad i\geq 2. \label{eq:telescope}
	\end{align}
 This follows from coupling two Markov chains $M^i_{k},M_{k}^{i-1}$ started at $i,i-1$ respectively, and sharing the same i.i.d update sequence $\lp X_{k}\rp$. The two processes are equal for all $k\geq 1$ if $X_{1}\geq i-1$ while $M^i_{1}=i-1,M^{i-1}_{1}=i-2$ if $X_{1}\leq i-2$, implying eq.  \eqref{eq:telescope}.  For all $i\geq 1$ we obtain
	\begin{align}
	e_{ij}&=e_{i-1,j}+\lp e_{1,j+i-1}-\phi\lp j+l-1\rp \rp \prod_{k=0}^{i-2}F\lp k\rp&
	\\&=\phi\lp j\rp +\sum_{l=1}^{i} \lp e_{1,j+l-1}-\phi\lp j+l-1\rp\rp \prod_{k=0}^{l -2}F\lp k\rp
	\\&=\phi\lp j\rp+ \sum_{l=1}^{i} \lp e_{1,j+l-1}-\phi\lp j+l-1\rp\rp \frac{\mu\lp 0\rp }{\mu\lf 0, l-1\rf }
	\end{align}
with convention $\prod_{i=0}^{-1}=1$. Inserting this into \eqref{eq:expectationformula} gives
	\begin{align}
		\qE_{\nu}\phi\lp \sigma\rp= S_{\nu}\lp 0\rp   \phi\lp 0\rp + \nu\lp 0\rp \qE_{0}\phi\lp \sigma\rp   + \sum_{i=1}^{\infty}\nu\lp i\rp  \sum_{l=1}^{i} \lp \qE_{1}\phi\lp \sigma+l-1\rp-\phi\lp l-1\rp\rp  \frac{\mu\lp 0\rp }{\mu\lf 0, l-1\rf }.
	\end{align}

	The result follows from adding $\pm \nu\lp 0\rp \phi\lp 0\rp $ and interchanging the two sums.
\end{p}
\begin{koro}
\label{s:markovkoro}
\textcolor{white}{123}\\[-35pt]
\begin{itemize}
	\item If $\qE X^{q+1}<\infty$ then an invariant distribution $\mu$ of $M$ exists and $\int x^{q} d\mu\lp x\rp<\infty.$  Also, if there is $c_{X}>0$ such that  $\exp\lp c_{X}X\rp<\infty$ then $\int \exp\lp cx\rp d\mu\lp x\rp $ for all $c<c_{X}$ and $M$ is geometrically ergodic.
\item
Assume still that $\qE X^{q+1}<\infty$.  It holds that $\qE_{\nu} \sigma^{q^*\wedge \lp q+1\rp } <\infty $ for all $q^*\geq 0$ and all initial measures $\nu$ with $q^{*}$'th moment. Moreover, if $\nu$ and $X$ have exponential moment then so has $\sigma$.
\end{itemize}
\end{koro}
\begin{p}
The 1st point follows directly from \cref{s:lyapounov}.

To prove the 2nd point notice that we can without loss of generalization assume $X,M_{0}$ has support on $\N_{0}$ by replacing $X$ and $M_{0}$ with $\lceil X \rceil$ and $\lc M_{0}\rc.$

 We start with the power-moment case. For $z\in \R$ write $e^z_{i,j}$ for the variable $e_{i,j}$ with $\phi\lp x\rp= x^z$.  Write $q=r+n$, $r\in \lf 0,1\rp$. We show by induction over  $m=0,\dots, n+1$ that $\qE_{\nu} \sigma^{q_{m}\wedge q^*}<\infty$ with $q_{m}=r+m$, and for all $q^*\geq 0$ and $\nu$  with $q^*$'th moment.
 \\\\
The induction start $n=0$; if $r=0$ or $q^*=0$ the claim is trivial. Otherwise, note that $e^r_{0,0}\leq e^1_{0,0}<\infty$ by Kac's theorem. We can apply \cref{s:hittingtimesum}  and the mean value theorem to obtain
\begin{align}\label{lolhat}
\qE_{\nu}\sigma^{ q^*\wedge r}\leq C + C\sum_{i = 0}^{\infty}i^{q^*\wedge r-1} S_{\nu}\lp i\rp<\infty
\end{align}
where $C>0$ is sufficiently large\\
Assume now that the induction claim holds for some $m\leq n$. Since $q_{m}\leq q$, and $\pi$ have $q$'th moment, we can use the induction assumption to see that $\qE_{\pi} \sigma^{q_{m}}<\infty$.  It is well known that $P_{\pi}\lp \sigma=i\rp=P_{0}\lp \sigma\geq i\rp$ for $i\geq 0$  (see section 10.3.1 \cite{MEYN}). It follows that  $e^{q_{m}+1}_{0,0}<\infty$ and hence also $e^{q_{m}+1}_{i,j}<\infty$. We may now apply  \cref{s:hittingtimesum} and the mean value theorem to obtain
\begin{align}\label{lolhat}
\qE_{\nu}\sigma^{\lp q_{m}+1\rp \wedge q^*}\leq C + C\sum_{i = 0}^{\infty}i^{\lp q_{m}+1\rp \wedge q^*-1} S_{\nu}\lp i\rp<\infty
\end{align}
where $C>0$ is sufficiently large. \\\\
Consider now the case where $\int \exp\lp c_{\nu}y\rp d\nu\lp y\rp <\infty$. For $c<c_{\nu}$ we consider the function $\phi\lp x\rp=\exp\lp cx\rp$.  Combining \cref{s:lyapounov} above and theorem 15.0.1 ii) in \cite{MEYN} we get that $\qE_{0}\phi\lp \inf\lb j>0: M_{j}\leq K\rb\rp<\infty$ for some $K>0$ large and small $c>0.$ It follows that $e_{0,0}<\infty$ for a possibly smaller $c$ and hence also $e_{i,j}<\infty$ for $i,j\in \N_{0}$. By dominated convergence, we can choose $C$ large so that
\begin{align}
\forall i\in \N_{0}: \qE_{1}\sigma\phi'\lp \sigma+i\rp \leq C \exp\lp c' i \rp
\end{align}
for all $c'>c.$ Choose now  $c' \in \lp c,c_{\nu}\rp$, insert the above inequality into \eqref{lolhat} and apply Markovs inequality to obtain the desired result.

.\end{p}

\section{Acknowledgment}
The work is part of the Dynamical Systems Interdisciplinary Network, University of Copenhagen. Author wish to thank PhD-supervisors Susanne Ditlevsen and Eva Löcherbach for great supervision, as well as fruitful and inspiring discussions about this article and probability theory in general.
\end{document}